\documentclass[10pt]{article}

\RequirePackage{amsmath}
\RequirePackage{amsfonts}
\RequirePackage{amssymb}
\RequirePackage{hyperref}
\RequirePackage{multirow}
\RequirePackage{algorithm}
\RequirePackage{algpseudocode}
\RequirePackage{graphicx}
\RequirePackage{xcolor}
\RequirePackage{amsthm}
\RequirePackage{amsmath}
\RequirePackage{amsfonts}
\RequirePackage{amssymb}
\RequirePackage{hyperref}
\RequirePackage{multirow}
\RequirePackage{algorithm}
\RequirePackage{algpseudocode}
\RequirePackage{graphicx}
\RequirePackage{xcolor}
\allowdisplaybreaks

\usepackage{dpr_fec,ifthen,subcaption,dsfont}

\usepackage{isetechreport}
\def\reportyear{19T}
\def\reportno{001}
\def\revisionno{2}
\def\originaldate{March 8, 2019}
\def\revisiondate{\today}
\coraltrue
\cvcrfalse

\begin{document}

\title{Limited-Memory BFGS with Displacement Aggregation\thanks{This material is based upon work supported by the National Science Foundation under grant numbers CCF--1618717 and CCF--1740796.}}

\author{Albert S.~Berahas\thanks{E-mail: albertberahas@gmail.com}}
\author{Frank E.~Curtis\thanks{E-mail: frank.e.curtis@gmail.com}}
\author{Baoyu Zhou\thanks{E-mail: baoyu.zhou@lehigh.edu}}
\affil{Department of Industrial and Systems Engineering, Lehigh University}
\titlepage

\maketitle

\begin{abstract}
  A displacement aggregation strategy is proposed for the curvature pairs stored in a limited-memory BFGS (a.k.a.~L-BFGS) method such that the resulting (inverse) Hessian approximations are equal to those that would be derived from a full-memory BFGS method.  This means that, if a sufficiently large number of pairs are stored, then an optimization algorithm employing the limited-memory method can achieve the same theoretical convergence properties as when full-memory (inverse) Hessian approximations are stored and employed, such as a local superlinear rate of convergence under assumptions that are common for attaining such guarantees.  To the best of our knowledge, this is the first work in which a local superlinear convergence rate guarantee is offered by a quasi-Newton scheme that does not either store all curvature pairs throughout the entire run of the optimization algorithm or store an explicit (inverse) Hessian approximation.  Numerical results are presented to show that displacement aggregation within an adaptive L-BFGS scheme can lead to better performance than standard L-BFGS.
\end{abstract}

\newcommand{\rev}[1]{#1}

\newcommand{\BFGS}{\textnormal{BFGS}}
\newcommand{\DFP}{\textnormal{DFP}}
\newcommand{\AggBFGS}{\texttt{Agg-BFGS}}

\section{Introduction}\label{sec.introduction}

Quasi-Newton methods---in which one computes search directions using (inverse) Hessian approximations that are set using iterate and gradient displacement information from one iteration to the next---represent some of the most effective algorithms for minimizing nonlinear objective functions.\footnote{Quasi-Newton methods offer the ability to update Hessian and/or inverse Hessian approximations, which is why we state \emph{inverse} parenthetically here.  For ease of exposition throughout the remainder of the paper, we often drop mention of the inverse, although in many cases it is the approximation of the inverse, not the Hessian approximation, that is used in practice.}  The main advantages of such methods can be understood by contrasting their computational costs, storage requirements, and convergence behavior with those of steepest descent and Newton methods.  Steepest descent methods require only first-order derivative (i.e., gradient) information, but only achieve a local linear rate of convergence.  Newton's method can achieve a faster (namely, quadratic) rate of local convergence, but at the added cost of forming and factoring second-order derivative (i.e., Hessian) matrices.  Quasi-Newton methods lie between the aforementioned methods; they only require gradient information, yet by updating and employing Hessian approximations they are able to achieve a local superlinear rate of convergence.  For many applications, the balance between cost/storage requirements and convergence rate offered by quasi-Newton methods makes them the most effective.

Davidon is credited as being the inventor of quasi-Newton methods \cite{Davi91}.  For further information on their theoretical properties, typically when coupled with line searches to ensure convergence, see, e.g., \cite{ByrdNoce89,ByrdNoceYuan87,DennMore74,DennSchn96,NoceWrig06,Pear69,Powe76,Ritt79,Ritt81}.  Within the class of quasi-Newton methods, algorithms that employ Broyden--Fletcher--Goldfarb--Shanno (BFGS) approximations of Hessian matrices have enjoyed particular success; see \cite{Broy70,Flet70,Gold70,Shan70}.  This is true when minimizing smooth objectives, as is the focus in the aforementioned references, but also when minimizing nonsmooth \cite{BonnGilbLemaSaga95,CurtQue13,CurtQue15,HaarMietMaek04,LewiOver13,MiffSunQi98,VlceLuks01} or stochastic  \cite{berahas2016multi,berahas2020robust,byrd2016stochastic,Curt16,gower2016stochastic,keskar2016adaqn,mokhtari2015global,schraudolph2007stochastic,wang2017stochastic} functions.

For solving large-scale problems, i.e., minimizing objective functions involving thousands or millions of variables, \emph{limited-memory} variants of quasi-Newton methods, such as the limited-memory variant of BFGS \cite{Noce80} (known as L-BFGS), have been successful for various applications.  The main benefit of a limited-memory scheme is that one need not store \rev{nor compute matrix-vector products with} the often-dense Hessian approximation; rather, one need only store a few pairs of vectors, known as \emph{curvature pairs}, with dimension equal to the number of variables\rev{, and one can compute matrix-vector products with the corresponding Hessian approximation with relatively low computational cost}.  That said, one disadvantage of contemporary limited-memory methods is that they do not enjoy the local superlinear convergence rate properties achieved by full-memory schemes.  Moreover, one does not know \emph{a priori} what number of pairs should be maintained to attain good performance when solving a particular problem.  Recent work \cite{BoggByrd19} has attempted to develop an adaptive limited-memory BFGS method in which the number of stored curvature pairs is updated dynamically within the algorithm.  However, while this approach has led to some improved empirical performance, it has not been coupled with any strengthened theoretical guarantees.

We are motivated by the following question: Is it possible to design an \emph{adaptive} limited-memory BFGS scheme that stores and employs only a \emph{moderate number} of curvature pairs, yet in certain situations can achieve the same theoretical convergence rate guarantees as a full-memory BFGS scheme (such as a local superlinear rate of convergence)?  We have not yet answered this question to the fullest extent possible.  That said, in this paper, we do answer the following question, which we believe is a significant first step: Is it possible to develop a \emph{limited-memory-type}\footnote{By \emph{limited-memory-type} BFGS algorithm, we mean one that stores and employs a finite set of curvature pairs rather than an explicit Hessian approximation.} BFGS scheme that behaves \emph{exactly} as a full-memory BFGS scheme---in the sense that the sequence of Hessian approximations from the full-memory scheme is represented through stored curvature pairs---while storing only a number of curvature pairs \emph{not exceeding the dimension of the problem}?  We answer this question in the affirmative by showing how one can use \emph{displacement aggregation} such that curvature pairs computed by the algorithm may be modified adaptively to capture the information that would be stored in a full-memory BFGS approximation.

A straightforward application of the strategy proposed in this paper might not always offer practical advantages over previously proposed full-memory or limited-memory schemes in all settings.  Indeed, to attain full theoretical benefits in some settings, one might have to store a number of curvature pairs up to the number of variables in the minimization problem, in which case the storage requirements and computational costs of our strategy might exceed those of a full-memory BFGS method.  That said, our displacement aggregation strategy and corresponding theoretical results establish a solid foundation upon which one may design practically efficient approaches.  We propose one such approach, leading to an adaptive L-BFGS scheme that we show to outperform a standard L-BFGS approach on a set of test problems.  For this and other such adaptive schemes, we argue that displacement aggregation allows one to maintain \emph{more history} in the same number of curvature pairs than in previously proposed limited-memory methods.

We conclude the paper with a discussion of how our approach can be extended to Davidon-Fletcher-Powell (DFP) quasi-Newton updating \cite{Davi91}, the challenges of extending it to the entire Broyden class of updates \cite{DennSchn96}, and further ways that displacement aggregation can be employed in practically efficient algorithms.

\subsection{Contributions}

We propose and analyze a \emph{displacement aggregation} strategy for modifying the displacement pairs stored in an L-BFGS scheme such that the resulting method behaves equivalently to full-memory BFGS.  In particular, we show that if a stored iterate displacement vector lies in the span of the other stored iterate displacement vectors, then the gradient displacement vectors can be modified in such a manner that the same Hessian approximation can be generated using the modified pairs with one pair (i.e, the one corresponding to the iterate displacement vector that is in the span of the others) being ignored.  In this manner, one can iteratively \emph{throw out} stored pairs and maintain at most a number of pairs equal to the dimension of the problem while generating the same sequence of Hessian approximations that would have been generated in a full-memory BFGS scheme.  Employed within a minimization algorithm, this leads to an L-BFGS scheme that has the same convergence behavior as full-memory BFGS; e.g., it can achieve a local superlinear rate of convergence.  To the best of our knowledge, this is the first work in which a local superlinear convergence rate guarantee is provided for a limited-memory-type quasi-Newton algorithm.  We also show how our techniques can be employed within an adaptive L-BFGS algorithm---storing and employing only a small number of pairs relative to the dimension of the problem---that can outperform standard L-BFGS.  We refer to our proposed ``aggregated BFGS'' scheme as \AggBFGS.

Our proposal of displacement aggregation should be contrasted with the idea of \emph{merging} information proposed in~\cite{KoldOLeaNaza98}.  In this work, the authors prove---when minimizing convex quadratics only---that one can maintain an algorithm that converges finitely if previous pairs are replaced by linear combinations of pairs.  This fact motivates a scheme proposed by the authors---when employing an L-BFGS-type scheme in general---of replacing two previous displacement pairs by a single ``sum'' displacement pair.  In their experiments, the authors' idea often leads to reduced linear algebra time, since fewer pairs are stored, but admittedly worse performance compared to standard L-BFGS in terms of function evaluations.  By contrast, our approach does not use predetermined weights to merge two pairs into one.  Our scheme aggregates information stored in any number of pairs using a strategy to ensure that no information is lost, even if the underlying function is not a convex quadratic.  Our experiments show that our approach can often result in reduced iteration and function evaluations compared to standard L-BFGS.

\subsection{Notation}

Let $\R{}$ denote the set of real numbers (i.e., scalars), let $\R{}_{\geq0}$ (resp.,~$\R{}_{>0}$) denote the set of nonnegative (resp.,~positive) real numbers, let $\R{n}$ denote the set of $n$-dimensional real vectors, and let $\R{m \times n}$ denote the set of $m$-by-$n$-dimensional real matrices.  Let $\mathbb{S}^n$ denote the symmetric elements of $\R{n \times n}$.  Let $\N{} := \{0,1,2,\dots\}$ denote the set of natural numbers.  Let $\|\cdot\| := \|\cdot\|_2$.

We motivate our proposed scheme in the context of solving
\bequation\label{prob.f}
  \min_{x\in\R{n}}\ f(x),
\eequation
where $f : \R{n} \to \R{}$ is continuously differentiable.  Corresponding to derivatives of~$f$, we define the gradient function $g : \R{n} \to \R{n}$ and Hessian function $H : \R{n} \to \mathbb{S}^n$.  In terms of an algorithm for solving~\eqref{prob.f}, we append a natural number as a subscript for a quantity to denote its value during an iteration of the algorithm; e.g., for each iteration number $k \in \N{}$, we denote $f_k := f(x_k)$.  That said, when discussing BFGS updating for the Hessian approximations employed within an algorithm, we often simplify notation by referring to a generic set of curvature pairs that may or may not come from consecutive iterations within the optimization algorithm.  In such settings, we clarify our notation at the start of each discussion.

\subsection{Organization}

In \S\ref{sec.background}, we provide background on BFGS updating.  In \S\ref{sec.aggregation}, we motivate, present, and analyze our displacement aggregation approach.  The results of numerical experiments with \AggBFGS{} are provided in~\S\ref{sec.numerical}.  Concluding remarks are given in \S\ref{sec.conclusion}.

\section{Background on BFGS}\label{sec.background}

The main idea of BFGS updating can be described as follows.  Let the $k$th iterate generated by an optimization algorithm  be denoted as $x_k$.  After an iterate displacement (or step) $s_k$ has been computed, one sets the subsequent iterate as $x_{k+1} \gets x_k + s_k$.  Then, in order to determine the subsequent step $s_{k+1}$, one uses the minimizer of the quadratic model $m_{k+1} : \R{n} \to \R{}$ given by
\bequationNN
  d_{k+1} \gets \arg\min_{d\in\R{n}} m_{k+1}(d),\ \ \text{where}\ \ m_{k+1}(d) = f_{k+1} + g_{k+1}^Td + \thalf d^T M_{k+1}d
\eequationNN
and $M_{k+1} \in \R{n\times n}$ is a Hessian approximation.  In a line search approach, for example, one computes $d_{k+1}$ by minimizing $m_{k+1}$, then computes $s_{k+1} \gets \alpha_{k+1}d_{k+1} = -\alpha_{k+1}M_{k+1}^{-1}g_{k+1}$, where $\alpha_{k+1}$ is chosen by a line search for $f$ from $x_{k+1}$ along $d_{k+1}$.  With $f_{k+1}$ and $g_{k+1}$ determined by evaluations at $x_{k+1}$, all that remains toward specifying the model $m_{k+1}$ is to choose $M_{k+1}$.  In a quasi-Newton method~\cite{Davi91}, one chooses $M_{k+1}$ such that it is symmetric---i.e., like the exact Hessian $H_{k+1}$, it is an element of $\mathbb{S}^n$---and satisfies the \emph{secant equation}
\bequation\label{eq.secant}
  M_{k+1}s_k = g_{k+1} - g_k =: y_k.
\eequation
Specifically in BFGS, one computes $W_{k+1} := M_{k+1}^{-1}$ to solve
\bequationNN
  \min_{W\in\mathbb{S}^n}\ \|W - W_k\|_\Mcal\ \ \st\ \ W = W^T\ \ \text{and}\ \ Wy_k = s_k,
\eequationNN
where $\|\cdot\|_\Mcal$ is a weighted Frobenius norm with weights defined by any matrix $\Mcal$ satisfying the secant equation (or, for concreteness, one can imagine the weight matrix being the \emph{average Hessian} between $x_k$ and $x_{k+1}$).  If one chooses $M_0 \succ 0$ and ensures that $s_k^Ty_k > 0$ for all $k \in \N{}$, as is often done by employing a (weak) Wolfe line search, then it follows that $M_k \succ 0$ for all $k \in \N{}$; see \cite{NoceWrig06}.

Henceforth, let us focus on the sequence of inverse Hessian approximations $\{W_k\}$.  After all, this sequence, not $\{M_k\}$, is the one that is often computed in practice since, for all $k \in \N{}$, the minimizer of $m_k$ can be computed as $d_k \gets -W_kg_k$.  Moreover, all of our discussions about the inverse Hessian approximations $\{W_k\}$ have corresponding meanings in terms of $\{M_k\}$ since $W_k \equiv M_k^{-1}$ for all $k \in \N{}$.

\subsection{Iterative and Compact Forms}

As is well known, there are multiple ways to construct or merely compute a matrix-vector product with a BFGS inverse Hessian approximation \cite{NoceWrig06}.  For our purposes, it will be convenient to refer to two ways of constructing such approximations: an iterative form and a compact form \cite{ByrdNoceSchn94}.  For convenience, let us temporarily drop from our notation the dependence on the iteration number of the optimization algorithm and instead talk generically about constructing an inverse Hessian approximation $\Wbar \succ 0$.  Regardless of whether one employs all displacements pairs since the start of the run of the optimization algorithm (leading to a \emph{full-memory} approximation) or only a few recent pairs (leading to a \emph{limited-memory} approximation), the approximation can be thought of as being constructed from some initial approximation $W \succ 0$ and a set of iterate and gradient displacements such that all iterate/gradient displacement inner products are positive, i.e.,
\bsubequations\label{eq.pairs}
  \begin{align}
    S &= \bbmatrix s_1 & \cdots & s_m \ebmatrix \in \R{n \times m} \\ \text{and}\ \ 
    Y &= \bbmatrix y_1 & \cdots & y_m \ebmatrix \in \R{n \times m}
  \end{align}
\esubequations
where
\bequation\label{eq.rho}
  \rho = \bbmatrix \frac{1}{s_1^Ty_1} & \cdots & \frac{1}{s_m^Ty_m} \ebmatrix^T \in \R{m}_{>0}.
\eequation

Algorithm~\ref{alg.lbfgs_iterative} computes the BFGS inverse Hessian approximation from an initial matrix $W \succ 0$ and the displacement pairs in \eqref{eq.pairs}.  In the context of an optimization algorithm using BFGS updating, this matrix would be set with a single update in each iteration rather than re-generated from scratch in every iteration using historical iterate/gradient displacements.  However, we write the strategy in this iterative form for convenience of our analysis.  The output $\Wbar$ represents the matrix to be employed in the step computation of the (outer) optimization algorithm.  As a function of the inputs, we denote the output as $\Wbar=\BFGS(W,S,Y)$.

\balgorithm[ht]
  \caption{:\ BFGS Matrix Construction, Iterative Form}
  \label{alg.lbfgs_iterative}
  \balgorithmic[1]
    \Require $W \succ 0$ and $(S,Y)$ as in \eqref{eq.pairs}, with $\rho$ as in \eqref{eq.rho}.
    \State Initialize $\Wbar \gets W$.
    \For{$j = 1,\dots,m$}
      \State Set
      \bsubequations\label{eq.W_update}
        \begin{align}
          U_j &\gets I - \rho_jy_js_j^T, \\
          V_j &\gets \rho_js_js_j^T, \\ \text{and}\ \ 
          \Wbar &\gets U_j^T\Wbar U_j + V_j.
        \end{align}
      \esubequations
    \EndFor
    \State \Return $\Wbar \equiv \textrm{BFGS}(W,S,Y)$
  \ealgorithmic
\ealgorithm

The updates performed in Algorithm~\ref{alg.lbfgs_iterative} correspond to a set of projections and corresponding corrections.  In particular, for each $j \in \{1,\dots,m\}$, the update \emph{projects out} curvature information along the step/direction represented by~$s_j$, then applies a subsequent \emph{correction} based on the gradient displacement represented by $y_j$; see \cite[Appendix~B]{CurtRobiZhou19}.  In this manner, one can understand the well-known fact that each update in a BFGS scheme involves a rank-two change of the matrix.

Rather than apply the updates iteratively, it has been shown that one can instead construct the BFGS approximation from the initial approximation by combining all low-rank changes directly.  The scheme in Algorithm~\ref{alg.lbfgs_compact}, which shows the compact form of the updates, generates the same output as Algorithm~\ref{alg.lbfgs_iterative}; see~\cite{ByrdNoceSchn94}.

\balgorithm[ht]
  \caption{:\ BFGS Matrix Construction, Compact Form}
  \label{alg.lbfgs_compact}
  \balgorithmic[1]
    \Require $W \succ 0$ and $(S,Y)$ as in \eqref{eq.pairs}, with $\rho$ as in \eqref{eq.rho}.
    \State Set $(R,D) \in \R{m \times m} \times \R{m \times m}$ with $R_{i,j} \gets s_i^Ty_j$ for all $(i,j)$ such that $1 \leq i \leq j \leq m$ and $D_{i,i} \gets s_i^Ty_i$ for all $i \in \{1,\dots,m\}$ (with all other elements being zero), i.e.,
    \bequation\label{eq.RD}
      R \gets \bbmatrix s_1^Ty_1 & \cdots & s_1^Ty_m \\ & \ddots & \vdots \\ & & s_m^Ty_m \ebmatrix\ \ \text{and}\ \ D \gets \bbmatrix s_1^Ty_1 & & \\ & \ddots & \\ & & s_m^Ty_m \ebmatrix.
    \eequation
    \State Set
    \bequation\label{eq.compact}
      \Wbar \gets W + \bbmatrix S & WY \ebmatrix \bbmatrix R^{-T}(D + Y^TWY)R^{-1} & -R^{-T} \\ -R^{-1} & 0 \ebmatrix \bbmatrix S^T \\ Y^TW \ebmatrix.
    \eequation
    \State \Return $\Wbar \equiv \BFGS(W,S,Y)$
  \ealgorithmic
\ealgorithm

We note in passing that for computing a matrix-vector product with an BFGS (or L-BFGS) approximation without constructing the approximation matrix itself, it is well-known that one can use the so-called two-loop recursion \cite{Noce80}.

\subsection{Convergence and Local Rate Guarantees}

An optimization algorithm using a BFGS scheme for updating Hessian approximations can be shown to converge and achieve a local superlinear rate of convergence under assumptions that are common for attaining such guarantees.  This theory for BFGS relies heavily on so-called \emph{bounded deterioration} of the updates or the \emph{self-correcting} behavior of the updates, and, for attaining a fast local rate of convergence, on the ability of the updating scheme to satisfy the well-known \emph{Dennis-Mor\'e condition} for superlinear convergence.  See~\cite{DennMore74,DennSchn96} for further information.

We show in the next section that our \AggBFGS{} approach generates the same sequence of matrices as full-memory BFGS.  Hence, an optimization algorithm that employs our updating scheme maintains \emph{all} of the convergence and convergence rate properties of full-memory BFGS.  For concreteness, we state one such result as the following theorem\rev{, which follows from Theorems~6.5 and 6.6} in \cite{NoceWrig06}.  We cite this result later in the paper to support our claim that our \AggBFGS{} scheme is a limited-memory-type quasi-Newton approach that can be used to achieve local superlinear convergence for an optimization algorithm.

\btheorem\label{th.bfgs_superlinear}
  Suppose that $f$ is twice continuously differentiable and that one employs an algorithm that generates a sequence of iterates $\{x_k\}$ according to
  \bequationNN
    d_k \gets -W_kg_k\ \ \text{and}\ \ x_{k+1} \gets x_k + \alpha_kd_k,
  \eequationNN
  where $\{W_k\}$ is generated using the BFGS updating scheme and, for all $k \in \N{}$, the stepsize $\alpha_k \in \R{}_{>0}$ is computed from a line search \rev{(initialized with a unit stepsize)} to satisfy the Armijo-Wolfe conditions; see equation~(3.6) in~\cite{NoceWrig06}.  In addition, suppose that the algorithm converges to a point $x_* \in \R{n}$ at which the Hessian is \rev{positive definite and} Lipschitz continuous.  Then, $\{x_k\}$ converges to $x_*$ at a superlinear rate.
\etheorem

Such a result cannot be proved for L-BFGS \cite{Noce80}.  Common theoretical results for L-BFGS merely show that if \rev{the pairs employed have $s_k^Ty_k$ sufficiently positive and bounded above for all $k \in \N{}$, then} the Hessian approximations are sufficiently positive definite and bounded \rev{and} one can achieve a local linear rate of convergence, i.e., no better than the rate offered by a steepest descent method.

\section{Displacement Aggregation}\label{sec.aggregation}

We begin this section by proving a simple, yet noteworthy result about a consequence that occurs when one makes consecutive BFGS updates with iterate displacements that are linearly dependent.  We use this result and other empirical observations to motivate our proposed approach, which is stated in this section.  We close this section by proving that our approach is well defined, and we discuss how to implement it in an efficient manner.

\subsection{Motivation: Parallel Consecutive Iterate Displacements}

The following theorem shows that if one finds in BFGS that an iterate displacement is a multiple of the previous one, then one can skip the update corresponding to the prior displacement and obtain the same inverse Hessian approximation.

\btheorem\label{th.1_update_skip}
  Let $(S,Y)$ be defined as in \eqref{eq.pairs}, with $\rho$ as in \eqref{eq.rho}, and suppose that $s_j = \tau s_{j+1} $ for some $j \in \{1,\dots,m-1\}$ and some nonzero $\tau \in \R{}$.  Then, with
  \bequationNN
    \baligned
      \Stilde &:= \bbmatrix s_1 & \cdots & s_{j-1} & s_{j+1} & \cdots s_m \ebmatrix \\ \text{and}\ \ 
      \Ytilde &:= \bbmatrix y_1 & \cdots & y_{j-1} & y_{j+1} & \cdots y_m \ebmatrix, \\
    \ealigned
  \eequationNN
  Algorithm~\ref{alg.lbfgs_iterative} $($and~\ref{alg.lbfgs_compact}$)$ yields $\BFGS(W,S,Y) = \BFGS(W,\Stilde,\Ytilde)$ for any $W \succ 0$.
\etheorem
\bproof
  Let $W \succ 0$ be chosen arbitrarily.  For any $j \in \{1,\dots,m-1\}$, let $W_{1:j} := \BFGS(W,S_{1:j},Y_{1:j})$ where $S_{1:j} := \bbmatrix s_1 & \cdots & s_j \ebmatrix$ and $Y_{1:j} := \bbmatrix y_1 & \cdots & y_j \ebmatrix$.  By \eqref{eq.W_update},
  \bequation\label{eq.W_j+1}
    W_{1:j+1} = U_{j+1}^TU_j^TW_{1:j-1}U_jU_{j+1} + U_{j+1}^TV_jU_{j+1} + V_{j+1}.
  \eequation
  Since $s_j = \tau s_{j+1}$, it follows that
  \bequationNN
    \baligned
      &\ U_jU_{j+1} \\
      =&\ (I - \rho_jy_js_j^T)(I - \rho_{j+1}y_{j+1}s_{j+1}^T) \\
      =&\ \(I - \(\frac{1}{\tau s_{j+1}^Ty_j}\)\tau y_j s_{j+1}^T\)(I - \rho_{j+1}y_{j+1}s_{j+1}^T) \\
      =&\ I - \(\frac{1}{s_{j+1}^Ty_j}\)y_js_{j+1}^T - \rho_{j+1}y_{j+1}s_{j+1}^T + \(\frac{1}{s_{j+1}^Ty_j}\)\rho_{j+1}y_j(s_{j+1}^Ty_{j+1})s_{j+1}^T \\
      =&\ I - \rho_{j+1}y_{j+1}s_{j+1}^T = U_{j+1}
    \ealigned
  \eequationNN
  and that
  \bequationNN
    \baligned
      V_jU_{j+1}
        &= (\rho_js_js_j^T)(I - \rho_{j+1}y_{j+1}s_{j+1}^T) \\
        &= \(\frac{1}{\tau s_{j+1}^Ty_j}\)\tau^2s_{j+1}s_{j+1}^T(I - \rho_{j+1}y_{j+1}s_{j+1}^T) \\
        &= \(\frac{1}{s_{j+1}^Ty_j}\)\tau \(s_{j+1}s_{j+1}^T - \rho_{j+1} s_{j+1}(s_{j+1}^Ty_{j+1})s_{j+1}^T\) = 0.
    \ealigned
  \eequationNN
  Combining these facts with \eqref{eq.W_j+1}, one finds that
  \bequationNN
    W_{1:j+1} = U_{j+1}^TW_{1:j-1}U_{j+1} + V_{j+1},
  \eequationNN
  meaning that one obtains the same matrix by applying the updates corresponding to $(s_j,y_j)$ and $(s_{j+1},y_{j+1})$ as when one skips the update for $(s_j,y_j)$ and only applies the one for $(s_{j+1},y_{j+1})$.  The result now follows by the fact that, if one starts with the same initial matrix $W_{1:j+1}$, then applying the updates defined by \eqref{eq.W_update} with the same pairs $\{(s_{j+2},y_{j+2}),\dots,(s_m,y_m)\}$ yields the same result. \qed
\eproof

Theorem~\ref{th.1_update_skip} is a consequence of the \emph{over-writing} process that signifies BFGS.  That is, with the update associated with each curvature pair, the BFGS update over-writes curvature information along the direction corresponding to each iterate displacement.  The theorem shows that if two consecutive iterate displacement vectors are linearly dependent, then the latter update completely over-writes the former---regardless of the gradient displacement vectors---and as a result the same matrix is derived if the former update is skipped.

\begin{figure}[ht]
  \centering
  \begin{subfigure}[b]{0.48\textwidth}
    \includegraphics[width=\textwidth]{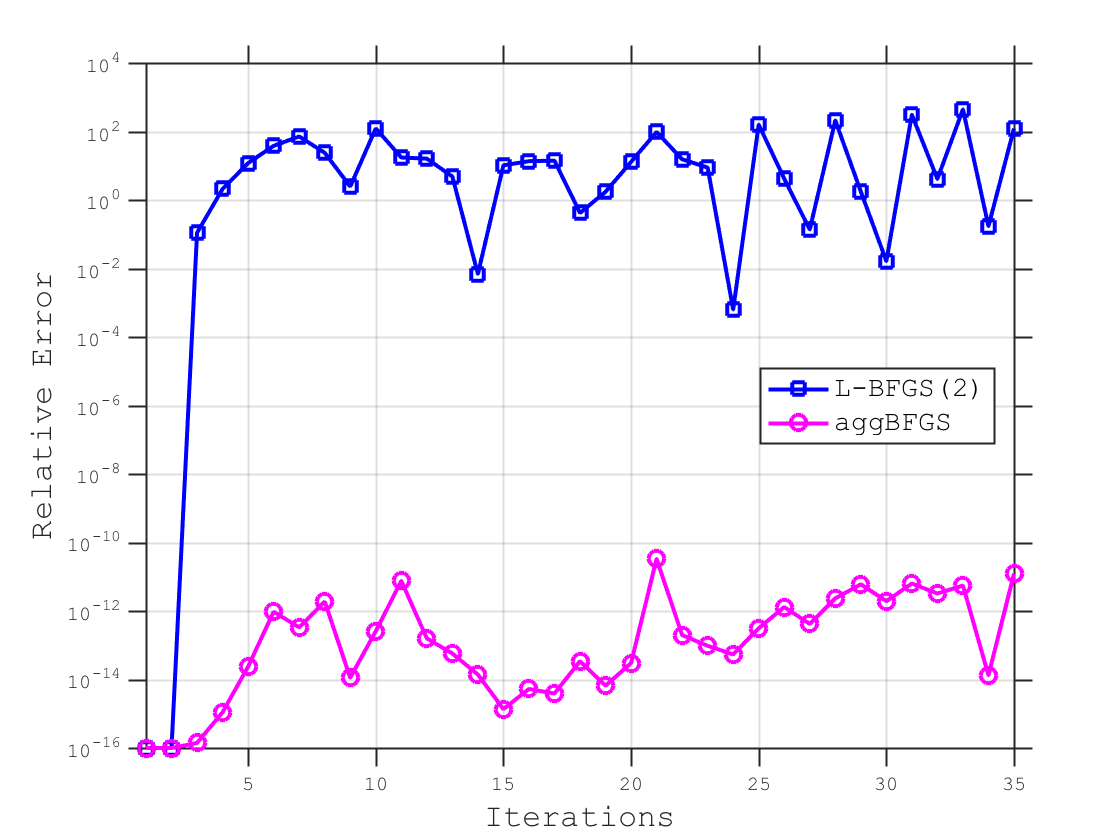}
    \caption{L-BFGS(2) and \AggBFGS{} vs.~BFGS}
    \label{fig.lbfgs_vs_aggbfgs}
  \end{subfigure}
  ~
  \begin{subfigure}[b]{0.48\textwidth}
    \includegraphics[width=\textwidth]{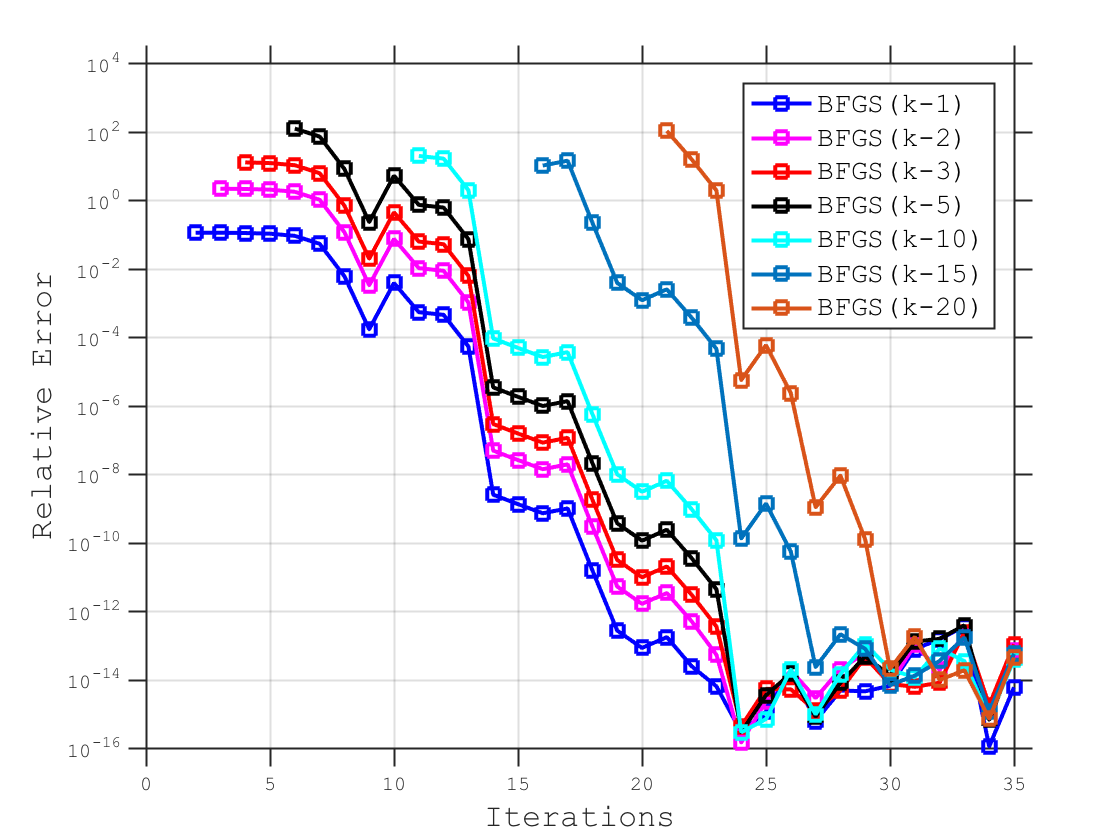}
    \caption{BFGS($k-j$) vs.~BFGS}
    \label{fig.bfgs_k}
  \end{subfigure}
  \caption{Illustrations that the over-writing of curvature information in BFGS inverse Hessian approximation updating is not absolute, even when previous iterate displacements lie in the span of subsequent iterate displacements. Relative error is the maximum absolute difference between corresponding matrix entries divided by the largest absolute value of an element of the BFGS inverse Hessian approximation.}
  \label{fig.bfgs_overwriting}
\end{figure}

What if a previous iterate displacement is not parallel with a subsequent displacement, but is in the span of multiple subsequent displacements?  One might suspect that the information in the previous displacement might get over-written and can be ignored.  \emph{It is not so simple}.  We illustrate this in two ways.
\bitemize
  \item Suppose that one accumulates curvature pairs $\{(s_j,y_j)\}_{j=0}^k$ and compares the corresponding BFGS approximations with those generated by an L-BFGS scheme with a history length of $n$.  (Hereafter, we refer to the latter scheme as L-BFGS($n$).)  Suppose also that for $k > n$ one finds that the latest $n$ iterate displacements span $\R{n}$.  If this fact meant that the updates corresponding to these pairs would over-write all curvature information from previous pairs, then one would find no difference between BFGS and L-BFGS($n$) approximations for $k > n$. However, one does not find this to be the case.  In Figure~\ref{fig.lbfgs_vs_aggbfgs}, we plot the maximum absolute difference between corresponding matrix entries of the BFGS and L-BFGS($2$) inverse Hessian approximations divided by the largest absolute value of an element of the BFGS inverse Hessian approximation when an algorithm is employed to minimize the Rosenbrock function \cite{Rose60}
  \bequation
    f(x_{(1)},x_{(2)}) = 100(x_{(2)} - x_{(1)}^2)^2 + (1 - x_{(1)})^2.
  \eequation
The pairs are generated by running an optimization algorithm using BFGS updating with an Armijo-Wolfe line search; the different approximations were computed as side computations.  It was verified that for $k > 2$ the iterate displacements used to generate the L-BFGS(2) matrices were linearly independent.  One finds that the differences between the approximations is large for $k > n$.  For the purposes of comparison, we also plot the differences between the BFGS inverse Hessian approximations and those generated by our \AggBFGS{} scheme described later; these are close to machine precision.
  \item Another manner in which one can see that historical information contained in a BFGS approximation is not always completely over-written once subsequent iterate displacements span $\R{n}$ is to compare BFGS approximations with those generated by a BFGS scheme that starts \emph{$j$ iterations late}.  Let us refer to such a scheme as  BFGS$(k-j)$.  For example, for $j=1$, BFGS$(k-1)$ approximations employ all pairs since the beginning of the run of the algorithm \emph{except} the first one.  By observing the magnitude of the differences in the approximations as~$k$ increases, one can see how long the information from the first pair \emph{lingers} in the BFGS approximation.  In Figure~\ref{fig.bfgs_k}, we plot the maximum absolute difference between corresponding matrix entries of the BFGS and BFGS$(k-j)$ inverse Hessian approximations divided by the largest absolute value of an element of the BFGS inverse Hessian approximation for various values of~$j$ using the same pairs generated by the algorithm from the previous bullet.  We only plot once the corresponding matrices start to deviate, i.e., for the BFGS and BFGS$(k-j)$ approximations the matrices only start to differ at the $j+1$st iteration. One finds, e.g., that the information from the first pair lingers---in the sense that it influences the BFGS approximation---for iteration numbers beyond $n=2$.
\eitemize

On the other hand, one might expect that the information stored in a full-memory BFGS matrix \emph{could be} contained in a set of at most $n$ curvature pairs, although not necessarily a subset of the pairs that are generated by the optimization algorithm.  For example, this might be expected by recalling the fact that if a BFGS method with an exact line search is used to minimize a strongly convex quadratic function, then the Hessian of the quadratic can be recovered \emph{exactly} if~$n$ iterations are performed; see \cite{DennSchn96,KoldOLeaNaza98,NoceWrig06}.  Letting $W$ represent a BFGS matrix---perhaps computed from more than $n$ pairs---this means that one could run an auxiliary BFGS method to minimize $x^TWx$ to \emph{re-generate} $W$ using at most $n$ pairs.  By noting the equivalence between the iterates generated by this auxiliary BFGS scheme and the conjugate gradient (CG) method, one finds that the modified iterate displacements would lie in a certain Krylov subspace determined by the initial point and the matrix $W$ \cite{DennSchn96,NoceWrig06}.

Practically speaking, it would not be computationally efficient to run an entire auxiliary BFGS scheme (for minimizing $W$) in order to capture $W$ using a smaller number of curvature pairs.    Moreover, we are interested in a scheme that can be used to reduce the number of curvature pairs that need to be stored even when an iterate displacement lies in the span of, say, only $m < n$ subsequent displacements.  Our \AggBFGS{} scheme is one efficient approach for achieving this goal.

\subsection{Basics of \AggBFGS}\label{eq.basics_aggBFGS}

The basic building block of our proposed scheme can be described as follows.  Suppose that, in addition to the iterate/gradient displacement information $(S,Y) \equiv ([s_1\ \cdots\ s_m],[y_1\ \cdots\ y_m]) =: (S_{1:m},Y_{1:m})$ as defined in~\eqref{eq.pairs} with inner products yielding~\eqref{eq.rho}, one also has a previous curvature pair
\bequation\label{eq.pair0}
  (s_0,y_0) \in \R{n} \times \R{n}\ \ \text{with}\ \ \rho_0 = 1/s_0^Ty_0 > 0
\eequation
such that
\bequation\label{eq.dependence}
  s_0 = S_{1:m}\tau\ \ \text{for some}\ \ \tau \in \R{m}.
\eequation
(As in the proof of Theorem~\ref{th.1_update_skip}, we have introduced the subscript ``$1:m$'' for $S$ and~$Y$ to indicate the dependence of these matrices on quantities indexed $1$ through~$m$.  We make use of similar notation throughout the rest of the paper.)  Given the linear dependence of the iterate displacement $s_0$ on the iterate displacements in $S_{1:m}$, our goal is to determine \emph{aggregated} gradient displacements
\bequation\label{eq.Ytilde}
  \Ytilde_{1:m} := \bbmatrix \ytilde_1 & \cdots & \ytilde_m \ebmatrix
\eequation
such that, for a given initial matrix $W \succ 0$, one finds
\bequation\label{eq.equal}
  \BFGS(W,S_{0:m},Y_{0:m}) = \BFGS(W,S_{1:m},\Ytilde_{1:m}).
\eequation
That is, our goal is to determine $\Ytilde_{1:m}$ such that the matrix $\BFGS(W,S_{0:m},Y_{0:m})$ is equivalently generated by ignoring $(s_0,y_0)$ and employing $(S_{1:m},\Ytilde_{1:m})$.

\bremark
  For simplicity, we are presenting the basics of \AggBFGS{} using indices for the iterate and gradient displacement vectors from 0 to $m$.  However, these need not correspond exactly to the iterates of the outer optimization algorithm.  In other words, our strategy and corresponding theoretical results apply equally for $(S,Y) \equiv ([s_{k_1}\ \cdots\ s_{k_m}],[y_{k_1}\ \cdots\ y_{k_m}])$ where $\{k_1,\dots,k_m\}$ are some $($not necessarily consecutive$)$ iteration numbers in the outer optimization algorithm.  We remark on this generality further when presenting the implementation of our scheme within a complete optimization algorithm; see~\S\ref{eq.implementation}.
\eremark

A basic view of our approach is stated as Algorithm~\ref{alg.aggregate_basic}, wherein we invoke the compact form of BFGS updates (recall \eqref{eq.compact}).  However, while Algorithm~\ref{alg.aggregate_basic} provides an understanding of the intent of our displacement aggregation scheme, it does not specify the manner in which $\Ytilde_{1:m}$ yielding \eqref{eq.equal} can be found.  It also does not address the fact that multiple such matrices might exist.  Toward a specific scheme, one can expand the compact forms of $\BFGS(W,S_{0:m},Y_{0:m})$ and $\BFGS(W,S_{1:m},\Ytilde_{1:m})$, specify that the aggregated displacements have the form
\bequation\label{eq.Ytilde_formula}
  \Ytilde_{1:m} = W^{-1}S_{1:m} \bbmatrix A & 0 \ebmatrix + y_0 \bbmatrix b \\ 0 \ebmatrix^T + Y_{1:m}
\eequation
for some $A \in \R{m \times (m-1)}$ and $b \in \R{m-1}$, and compare like terms in order to derive three key equations that ensure that \eqref{eq.equal} holds.  The zero blocks in the variable matrices in~\eqref{eq.Ytilde_formula} are motivated by Theorem~\ref{th.1_update_skip}, since in the case of $m=1$ one can set $\Ytilde_{1:1} = Y_{1:1}$; otherwise, $A$ and $b$ need to be computed to satisfy \eqref{eq.equal}.

\balgorithm[ht]
  \caption{:\ Displacement Aggregation, Basic View}
  \label{alg.aggregate_basic}
  \balgorithmic[1]
    \Require $W \succ 0$, $(S_{1:m},Y_{1:m},\rho_{1:m})$ as in \eqref{eq.pairs}--\eqref{eq.rho}, and $(s_0,y_0,\rho_0)$ as in \eqref{eq.pair0}--\eqref{eq.dependence}.
    \State Set $\Ytilde_{1:m}$ as in \eqref{eq.Ytilde} such that, with
    \bequation\label{eq.RDtilde_basic}
      \Rtilde_{1:m} := \bbmatrix s_1^T\ytilde_1 & \cdots & s_1^T\ytilde_m \\ & \ddots & \vdots \\ & & s_m^T\ytilde_m \ebmatrix\ \ \text{and}\ \ \Dtilde_{1:m} := \bbmatrix s_1^T\ytilde_1 & & \\ & \ddots & \\ & & s_m^T\ytilde_m \ebmatrix,
    \eequation
    one finds with
    \bequationNN
      R_{0:m} := \bbmatrix s_0^Ty_0 & s_0^TY_{1:m} \\ & R_{1:m} \ebmatrix\ \ \text{and}\ \ D_{0:m} := \bbmatrix s_0^Ty_0 & \\ & D_{1:m} \ebmatrix
    \eequationNN
    that
    \bequationNN
      \baligned
          &\ \BFGS(W,S_{0:m},Y_{0:m}) \\
        \equiv &\ W + \bbmatrix S_{0:m} & WY_{0:m} \ebmatrix \bbmatrix R_{0:m}^{-T}(D_{0:m} + Y_{0:m}^TWY_{0:m})R_{0:m}^{-1} & -R_{0:m}^{-T} \\ -R_{0:m}^{-1} & 0 \ebmatrix \bbmatrix S_{0:m}^T \\ Y_{0:m}^TW \ebmatrix \\
        = &\ W + \bbmatrix S_{1:m} & W\Ytilde_{1:m} \ebmatrix \bbmatrix \Rtilde_{1:m}^{-T}(\Dtilde_{1:m} + \Ytilde_{1:m}^TW\Ytilde_{1:m})\Rtilde_{1:m}^{-1} & -\Rtilde_{1:m}^{-T} \\ -\Rtilde_{1:m}^{-1} & 0 \ebmatrix \bbmatrix S_{1:m}^T \\ \Ytilde_{1:m}^TW \ebmatrix \\
        \equiv &\ \BFGS(W,S_{1:m},\Ytilde_{1:m}).
      \ealigned
    \eequationNN
    \State \Return $\Ytilde_{1:m}$.
  \ealgorithmic
\ealgorithm

Following these steps, one obtains the more detailed view of our scheme that is stated in Algorithm~\ref{alg.aggregate_detailed}, where the three key equations are stated as \eqref{eq.key}.  Algorithm~\ref{alg.aggregate_detailed} does not specify a precise scheme for computing $(A,b)$ in order to satisfy~\eqref{eq.key}.  We specify such a scheme in \S\ref{eq.implementation} after first showing in the following subsection that real values for $(A,b)$ always exist to satisfy \eqref{eq.key}.

\balgorithm[ht]
  \caption{:\ Displacement Aggregation, Detailed View}
  \label{alg.aggregate_detailed}
  \balgorithmic[1]
    \Require $W \succ 0$, $(S_{1:m},Y_{1:m},\rho_{1:m})$ as in \eqref{eq.pairs}--\eqref{eq.rho}, and $(s_0,y_0,\rho_0)$ as in \eqref{eq.pair0}--\eqref{eq.dependence}.
    \State Set $\Ytilde_{1:m}$ as in \eqref{eq.Ytilde_formula} such that, with $\chi_0 := 1 + \rho_0\|y_0\|_W^2$, one finds
    \bsubequations\label{eq.key}
      \begin{align}
        \Rtilde_{1:m} =&\ R_{1:m}, \label{eq.R=Rtilde} \\
        \bbmatrix b \\ 0 \ebmatrix =&\ -\rho_0 (S_{1:m}^TY_{1:m} - R_{1:m})^T \tau,\ \ \text{and} \label{eq.b} \\
        (\Ytilde_{1:m} - Y_{1:m})^TW(\Ytilde_{1:m} - Y_{1:m}) =&\ \frac{\chi_0}{\rho_0} \bbmatrix b \\ 0 \ebmatrix \bbmatrix b \\ 0 \ebmatrix^T \label{eq.YWYtilde} \\
         &\ - \bbmatrix A & 0 \ebmatrix^T (S_{1:m}^TY_{1:m} - R_{1:m}) \nonumber \\
         &\ - (S_{1:m}^TY_{1:m} - R_{1:m})^T\bbmatrix A & 0 \ebmatrix. \nonumber
      \end{align}
    \esubequations
    \State \Return $\Ytilde_{1:m}$.
  \ealgorithmic
\ealgorithm

\subsection{Existence of Real Solutions for \AggBFGS{}}\label{sec.existence}

Our goal now is to prove the following theorem for our \AggBFGS{} scheme.

\btheorem\label{th.real}
  Suppose one has $W \succ 0$ along with:
  \bitemize
    \item $(S_{1:m},Y_{1:m})$ as defined in \eqref{eq.pairs} with the columns of $S_{1:m}$ being linearly independent and the vector $\rho_{1:m}$ as defined in \eqref{eq.rho}, and
    \item $(s_0,y_0,\rho_0)$ defined as in \eqref{eq.pair0} such that \eqref{eq.dependence} holds for some $\tau \in \R{m}$.
  \eitemize
  Then, there exist $A \in \R{m \times (m-1)}$ and $b \in \R{m-1}$ such that, with $\Ytilde_{1:m} \in \R{n \times m}$ defined as in \eqref{eq.Ytilde_formula}, the equations~\eqref{eq.key} hold.  Consequently, for this $\Ytilde_{1:m}$, one finds that
  \bequation\label{eq.products_remain_positive}
    s_i^T\ytilde_i = s_i^Ty_i > 0\ \ \text{for all}\ \ i \in \{1,\dots,m\}
  \eequation
  and $\BFGS(W,S_{0:m},Y_{0:m}) = \BFGS(W,S_{1:m},\Ytilde_{1:m})$.
\etheorem
\begin{proof}
First, observe that there are $(m+1)(m-1)$ unknowns in the formula~\eqref{eq.Ytilde_formula} for~$\Ytilde_{1:m}$; in particular, there are $m(m-1)$ unknowns in $A$ and $m-1$ unknowns in $b$, yielding $m(m-1) + (m-1) = (m+1)(m-1)$ unknowns in total.  To see the number of equations effectively imposed by \eqref{eq.key}, first notice that \eqref{eq.Ytilde_formula} imposes $\ytilde_m = y_m$.  Hence, by ignoring the last column of $R_{1:m}$ to define the submatrix
\bequation\label{eq.P}
  P := \bbmatrix s_1^Ty_1 & \cdots & s_1^Ty_{m-1} \\ 0 & \ddots & \vdots \\ \vdots & \ddots & s_{m-1}^Ty_{m-1} \\ 0 & \cdots & 0 \ebmatrix \in \R{m \times (m-1)},
\eequation
and similarly defining $\Ptilde$ with size $m \times (m-1)$ as a submatrix of $\Rtilde_{1:m}$, the key equations in~\eqref{eq.key} can be simplified to have the following form:
\bsubequations\label{eq.key_simplified}
  \begin{align}
    &\Ptilde =\ P, \label{eq.R=Rtilde_simplified} \\
    &b \ =\ -\rho_0 (S_{1:m}^TY_{1:m-1} - P)^T \tau, \label{eq.b_simplified} \\
    \text{and }\ \ \ & (\Ytilde_{1:m-1} - Y_{1:m-1})^TW(\Ytilde_{1:m-1} - Y_{1:m-1}) \nonumber \\
      &\qquad =\ \frac{\chi_0}{\rho_0} bb^T - A^T (S_{1:m}^TY_{1:m-1} - P) - (S_{1:m}^TY_{1:m-1} - P)^TA. \label{eq.YWYtilde_simplified}
  \end{align}
\esubequations
Observing the number of nonzero entries in \eqref{eq.R=Rtilde_simplified} (recall \eqref{eq.P}) and recognizing the symmetry in \eqref{eq.YWYtilde_simplified}, one finds that the effective number of equations are:
\bequationNN
  \baligned
      &\ \underbrace{m(m-1)/2}_\text{(in \eqref{eq.R=Rtilde_simplified})} +\underbrace{(m-1)}_\text{(in \eqref{eq.b_simplified})} +  \underbrace{m(m-1)/2}_\text{(in \eqref{eq.YWYtilde_simplified})} =  \underbrace{(m+1)(m-1)}_\text{(in \eqref{eq.key_simplified})}.
  \ealigned
\eequationNN
Hence, \eqref{eq.key_simplified} (and so \eqref{eq.key}) is a square system of linear and quadratic equations to be solved for the unknowns in the matrix $A$ and vector $b$.

Equation~\eqref{eq.b_simplified} represents a formula for $b \in \R{m-1}$.  Henceforth, let us assume that $b$ is equal to the right-hand side of this equation, meaning that all that remains is to determine that a real solution for $A$ exists.  Let us write
\bequationNN
  A = \bbmatrix a_1 & \cdots & a_{m-1} \ebmatrix\ \ \text{where $a_i$ has length $m$ for all}\ \ i \in \{1,\dots,m-1\}.
\eequationNN
Using this notation, \rev{one finds from \eqref{eq.Ytilde_formula} that the $j$th column of \eqref{eq.R=Rtilde_simplified} requires}
\bequationNN
  \rev{S_{1:j}^Ty_j = S_{1:j}^T\ytilde_j \iff S_{1:j}^Ty_j = S_{1:j}^T(W^{-1}S_{1:m}a_j + b_jy_0 + y_j);}
\eequationNN
\rev{hence, \eqref{eq.R=Rtilde_simplified} reduces to the system of affine equations}
\bequation\label{eq.R=Rtilde_j}
  S_{1:j}^TW^{-1}S_{1:m}a_j = - b_j S_{1:j}^T y_0\ \ \text{for all}\ \ j \in \{1,\dots,m-1\}.
\eequation
For each $j \in \{1,\dots,m-1\}$, let us write
\bequation\label{eq.aj}
  a_j = Q^{-1}\bbmatrix a_{j,1} \\ a_{j,2} \ebmatrix,\ \text{where}\ Q := S_{1:m}^TW^{-1}S_{1:m} \succ 0,
\eequation
with $a_{j,1}$ having length $j$ and $a_{j,2}$ having length $m-j$.  (Notice that $Q$ is positive definite under the conditions of the theorem, which requires that $S_{1:m}$ has full column rank.)  Then, in order for \eqref{eq.R=Rtilde_j} to be satisfied, one must have
\bequation\label{eq.aj1}
  a_{j,1} = -b_jS_{1:j}^Ty_0 \in \R{j}.
\eequation
Moreover, with this value for $a_{j,1}$, it follows that \eqref{eq.R=Rtilde_j}, and hence \eqref{eq.R=Rtilde_simplified}, is satisfied for any $a_{j,2}$.  Going forward, our goal is to show the existence of $a_{j,2} \in \R{m-j}$ for all $j \in \{1,\dots,m-1\}$ such that \eqref{eq.YWYtilde_simplified} holds, completing the proof.

Observing \eqref{eq.YWYtilde_simplified}, one finds with \eqref{eq.Ytilde_formula} that it may be written as
\bequation\label{eq.YWYtilde_compact}
  A^TQA + \Omega^TA + A^T\Omega - \omega\omega^T = 0,
\eequation
where
\bsubequations\label{eq.defs}
  \begin{align}
    \Omega &:= S_{1:m}^Ty_0b^T + S_{1:m}^TY_{1:m-1} - P \in \R{m \times (m-1)}, \\
    \text{and}\ \ \omega &:= \frac{b}{\sqrt{\rho_0}} \in \R{m-1}.
  \end{align}
\esubequations
One may rewrite equation~\eqref{eq.YWYtilde_compact} as
\bequation\label{eq.YWYtilde_compact_reduced}
  (QA + \Omega)^TQ^{-1}(QA + \Omega) = \omega\omega^T + \Omega^TQ^{-1}\Omega.
\eequation
Let us now rewrite the equations in \eqref{eq.YWYtilde_compact_reduced} in a particular form that will be useful for the purposes of our proof going forward.  Consider the matrix $QA + \Omega$ in \eqref{eq.YWYtilde_compact_reduced}.  By the definitions of $a_{j,1}$, $a_{j,2}$, $Q$, and $\Omega$ in \eqref{eq.aj}--\eqref{eq.defs}, as well as of $P$ from \eqref{eq.P}, the $j$th column of this matrix is given by
  \bequationNN
    \baligned
      \left[QA + \Omega\right]_j
        &= \bbmatrix -b_jS_{1:j}^Ty_0 \\ a_{j,2} \ebmatrix + \bbmatrix b_jS_{1:j}^Ty_0 \\ b_jS_{j+1:m}^Ty_0 \ebmatrix + \bbmatrix 0_j \\ S_{j+1:m}^Ty_j \ebmatrix \\
        &= \bbmatrix 0_j \\ a_{j,2} \ebmatrix + \bbmatrix 0_j \\ S_{j+1:m}^T(b_jy_0 + y_j) \ebmatrix,
    \ealigned
  \eequationNN
  where $0_j$ is a vector of zeros of length $j$.  Letting $L \in \R{m \times m}$ be any matrix such that $L^TL = Q^{-1}$ (whose existence follows since $Q^{-1} \succ 0$), defining $Z := [z_1\ \cdots\ z_{m-1}] \in \R{(m-1) \times (m-1)}$ as any matrix such that $Z^TZ = \omega\omega^T + \Omega^TQ^{-1}\Omega$ (whose existence follows since $\omega\omega^T + \Omega^TQ^{-1}\Omega \succeq 0$), and defining, for all $j \in \{1,\dots,m-1\}$,
  \bequation\label{def.phi}
    \phi_j(a_{j,2}) := L\(\bbmatrix 0_j \\ a_{j,2} \ebmatrix + \bbmatrix 0_j \\ S_{j+1:m}^T(b_jy_0 + y_j) \ebmatrix\),
  \eequation
  it follows that the $(i,j) \in \{1,\dots,m-1\} \times \{1,\dots,m-1\}$ element of \eqref{eq.YWYtilde_compact_reduced} is
  \bequation\label{eq.phi}
    \phi_i(a_{i,2})^T\phi_j(a_{j,2}) = z_i^Tz_j.
  \eequation
  
  Using the notation from the preceding paragraph, let us use an inductive argument to prove the existence of a real solution of \eqref{eq.YWYtilde_compact_reduced}.  This induction will follow the indices $\{1,\dots,m-1\}$ \emph{in reverse order}.  As a base case, consider the index $m-1$, in which case one has the one-dimensional unknown $a_{m-1,2}$.  One finds with
  \bequationNN
    a_{m-1,2}^* := -s_m^T(b_{m-1}y_0 + y_{m-1}) \in \R{}
  \eequationNN
  that
  \bequationNN
    \phi_{m-1}(a_{m-1,2}^*) = L\bbmatrix 0_{m-1} \\ -s_m^T(b_{m-1}y_0 + y_{m-1}) + s_m^T(b_{m-1}y_0 + y_{m-1}) \ebmatrix = 0_m.
  \eequationNN
  Hence, letting $a_{m-1,2} = a_{m-1,2}^* + \lambda_{m-1}$, where $\lambda_{m-1}$ is one-dimensional, one finds that the left-hand side of the $(i,j)=(m-1,m-1)$ equation in \eqref{eq.phi} is $\|\phi_{m-1}(a_{m-1,2})\|_2^2$.  This is a strongly convex quadratic in the unknown~$\lambda_{m-1}$.  Since $\|\phi_{m-1}(a_{m-1,2}^*)\|_2 = 0$ and $z_{m-1}^Tz_{m-1} \geq 0$, there exists~$\lambda_{m-1}^* \in \R{}$ such that $a_{m-1,2} = a_{m-1,2}^* + \lambda_{m-1}^* \in \R{}$ satisfies the $(i,j)=(m-1,m-1)$ equation in \eqref{eq.phi}.
  
  Now suppose that there exists real $\{a_{\ell + 1,2},\dots,a_{m-1,2}\}$ such that \eqref{eq.phi} holds for all $(i,j)$ with $i \in \{\ell+1,\dots,m-1\}$ and $j \in \{i,\dots,m-1\}$.  (By symmetry, these values also satisfy the $(j,i)$ elements of \eqref{eq.phi} for these same values of the indices $i$ and $j$.)  To complete the inductive argument, we need to show that this implies the existence of $a_{\ell,2} \in \R{m-\ell}$ satisfying \eqref{eq.phi} for all $(i,j)$ with $i \in \{j,\dots,m-1\}$ and $j = \ell$, i.e., solving the following system for $a_{\ell,2}$:
  \bsubequations\label{eq.linear_quadratic}
    \begin{align}
      \phi_{m-1}(a_{m-1,2})^T\phi_{\ell}(a_{\ell,2}) &= z_{m-1}^Tz_{\ell} \label{eq.affine_begin} \\
      &\ \vdots \nonumber \\
      \phi_{\ell+1}(a_{\ell+1,2})^T\phi_{\ell}(a_{\ell,2}) &= z_{\ell+1}^Tz_{\ell} \label{eq.affine_end} \\
      \phi_{\ell}(a_{\ell,2})^T\phi_{\ell}(a_{\ell,2}) &= z_{\ell}^Tz_{\ell}. \label{eq.quadratic}
    \end{align}
  \esubequations
  Notice that \eqref{eq.affine_begin}--\eqref{eq.affine_end} are affine equations in $a_{\ell,2}$, whereas \eqref{eq.quadratic} is a quadratic equation in $a_{\ell,2}$.  For all $t \in \{\ell+1,\dots,m-1\}$, let
  \bequationNN
    \psi_{\ell+1,t} := \bbmatrix 0_{t-(\ell+1)} \\ a_{t,2} + S_{t+1:m}^T(b_ty_0 + y_t) \ebmatrix \in \R{m-(\ell+1)},
  \eequationNN
  so that, by \eqref{def.phi}, one may write
  \bequationNN
    \phi_t(a_{t,2}) = L\bbmatrix 0_{\ell+1} \\ \psi_{\ell+1,t} \ebmatrix\ \ \text{for all}\ \ t \in \{\ell+1,\dots,m-1\}.
  \eequationNN
  Our strategy is first to find $a_{\ell,2}^* \in \R{m-\ell}$ satisfying \eqref{eq.affine_begin}--\eqref{eq.affine_end} such that
  \bequation\label{spanspace}
    a_{\ell,2}^* + S_{\ell+1:m}^T(b_{\ell}y_0 + y_{\ell}) \in \linspan\left \{\bbmatrix 0 \\ \psi_{\ell+1,\ell+1} \ebmatrix,\dots,\bbmatrix 0 \\ \psi_{\ell+1,m-1} \ebmatrix\right\}
  \eequation
  and (cf.~\eqref{eq.quadratic})
  \bequation\label{starting_point}
  \phi_{\ell}(a_{\ell,2}^*)^T\phi_{\ell}(a_{\ell,2}^*) \leq z_{\ell}^Tz_{\ell}.
  \eequation
  Once this is done, we will argue the existence of a nonzero vector $\abar_{\ell,2} \in \R{m-\ell}$ such that $a_{\ell,2}^* + \lambda_{\ell}\abar_{\ell,2}$ satisfies \eqref{eq.affine_begin}--\eqref{eq.affine_end} for arbitrary one-dimensional $\lambda_{\ell}$.  From here, it will follow by the fact that the left-hand side of \eqref{eq.quadratic} is a strongly convex quadratic in the unknown~$\lambda_{\ell}$ and the fact that \eqref{starting_point} holds that we can claim that there exists $\lambda_{\ell}^* \in \R{}$ such that $a_{\ell,2} = a_{\ell,2}^* + \lambda_{\ell}^* \abar_{\ell,2} \in \R{m-\ell}$ satisfies \eqref{eq.linear_quadratic}.
  
  To achieve the goals of the previous paragraph, first let $c$ be the column rank of $[\psi_{\ell+1,\ell+1}\ \cdots\ \psi_{\ell+1,m-1}]$ so that there exists $\{t_1,\dots,t_c\} \subseteq \{\ell+1,\dots,m-1\}$ with
  \bequation\label{spanspace_new}
    \linspan\left \{\bbmatrix 0 \\ \psi_{\ell+1,t_1} \ebmatrix,\dots,\bbmatrix 0 \\ \psi_{\ell+1,t_c} \ebmatrix\right\} = \linspan\left \{\bbmatrix 0 \\ \psi_{\ell+1,\ell+1} \ebmatrix,\dots,\bbmatrix 0 \\ \psi_{\ell+1,m-1} \ebmatrix\right\}.
  \eequation
  For completeness, let us first consider the extreme case when $c=0$.  In this case,
  \bequation\label{psi.0}
    \psi_{\ell+1,t} = 0_{m-(\ell+1)}\ \ \text{and}\ \ \phi_t(a_{t,2}) = 0_m\ \ \text{for all}\ \ t\in\{\ell+1,\ldots,m-1\}.
  \eequation
  Hence, by our induction hypothesis, it follows from \eqref{eq.phi} and \eqref{psi.0} that
  \bequation\label{z_t.0}
    z_t = 0_{m-1}\ \ \text{for all}\ \ t\in\{\ell+1,\ldots,m-1\}.
  \eequation
  Consequently from \eqref{psi.0} and \eqref{z_t.0}, the affine equations \eqref{eq.affine_begin}--\eqref{eq.affine_end} are satisfied by any $a_{\ell,2}^*\in\R{m-\ell}$.  In particular, one can choose 
  \bequationNN
    a_{\ell,2}^* = -S_{\ell+1:m}^T(b_{\ell}y_0 + y_{\ell}),
  \eequationNN
  and find by \eqref{def.phi} that
  \bequation
    \phi_{\ell}(a_{\ell,2}^*) = L\(\bbmatrix 0_{\ell} \\ a_{\ell,2}^* + S_{\ell+1:m}^T(b_{\ell}y_0 + y_{\ell}) \ebmatrix\) = 0_m,
  \eequation
  which shows that this choice satisfies \eqref{starting_point}.  Now consider the case when $c>0$. For $a_{\ell,2}^*$ to satisfy \eqref{spanspace}, it follows with \eqref{spanspace_new} that we must have
  \bequation\label{beta}
    a_{\ell,2}^* + S_{\ell+1:m}^T(b_{\ell}y_0 + y_{\ell}) = \bbmatrix 0 & \cdots & 0 \\ \psi_{\ell+1,t_1} & \cdots & \psi_{\ell+1,t_c} \ebmatrix\beta_{\ell},
  \eequation
  where $\beta_{\ell}$ has length $c$.  Choosing 
  \bequation\label{beta_value}
    \beta_{\ell} := \(\bbmatrix 0_{\ell+1} & \cdots & 0_{\ell+1} \\ \psi_{\ell+1,t_1} & \cdots & \psi_{\ell+1,t_c} \ebmatrix^TL^TL\bbmatrix 0_{\ell+1} & \cdots & 0_{\ell+1} \\ \psi_{\ell+1,t_1} & \cdots & \psi_{\ell+1,t_c} \ebmatrix\)^{-1} \bbmatrix z_{t_1}^T \\ \vdots \\ z_{t_c}^T \ebmatrix z_{\ell} \in \R{c},
  \eequation
  it follows with \eqref{beta} that, for any $t \in \{t_1,\dots,t_c\}$, one finds
  \bequation\label{eq.affine_1}
  \baligned
  \phi_{t}(a_{t,2})^T\phi_{\ell}(a_{\ell,2}^*) &= \bbmatrix 0_{\ell+1} \\ \psi_{\ell+1,t} \ebmatrix^TL^TL\bbmatrix 0_{\ell} \\ a_{\ell,2}^* + S_{\ell+1:m}^T(b_{\ell}y_0 + y_{\ell}) \ebmatrix \\
  &= \bbmatrix 0_{\ell+1} \\ \psi_{\ell+1,t} \ebmatrix^TL^TL\bbmatrix 0_{\ell+1} & \cdots & 0_{\ell+1} \\ \psi_{\ell+1,t_1} & \cdots & \psi_{\ell+1,t_c} \ebmatrix \beta_{\ell} = z_t^Tz_{\ell}.
  \ealigned
  \eequation
  We shall now prove that, for any $t \in \{\ell+1,\dots,m-1\}\backslash\{t_1,\dots,t_c\}$, one similarly finds that $\phi_t(a_{t,2})^T\phi_{\ell}(a_{\ell,2}^*) = z_{t}^Tz_{\ell}$.  Toward this end, first notice that for any such~$t$ it follows from \eqref{spanspace_new} that $\psi_{\ell+1,t} = [\psi_{\ell+1,t_1}\ \cdots\ \psi_{\ell+1,t_c}] \gamma_{\ell,t}$ for some $\gamma_{\ell,t} \in \R{c}$.  Combining the relationship~\eqref{spanspace_new} along with the inductive hypothesis that, for any pair $(i,j)$ with $i\in\{\ell+1,\dots,m-1\}$ and $j\in\{i,\dots,m-1\}$, one has
  \bequation\label{main_assumption}
  \bbmatrix 0_{\ell+1} \\ \psi_{\ell+1,i} \ebmatrix^TL^TL\bbmatrix 0_{\ell+1} \\ \psi_{\ell+1,j} \ebmatrix = \phi_i(a_{i,2})^T\phi_j(a_{j,2}) = z_i^Tz_j,
  \eequation
  it follows with the positive definiteness of $Q^{-1} = L^TL$ that
  \bequation\label{rank}
  \baligned
      \rank\(\bbmatrix z_{\ell+1} & \cdots & z_{m-1} \ebmatrix\) = &\ \rank\(\bbmatrix \phi_{\ell+1}(a_{\ell+1,2}) & \cdots & \phi_{m-1}(a_{m-1,2}) \ebmatrix\)\\
    =&\ \rank\(\bbmatrix \phi_{t_1}(a_{t_1,2}) & \cdots & \phi_{t_c}(a_{t_c,2}) \ebmatrix\)\\
    =&\ \rank\(\bbmatrix z_{t_1} & \cdots & z_{t_c} \ebmatrix\) = c.
  \ealigned
  \eequation
  From \eqref{rank}, it follows that for any $t \in \{\ell+1,\dots,m-1\}\backslash\{t_1,\dots,t_c\}$ there exists $\bar{\gamma}_{\ell,t}\in\R{c}$ such that $z_t = [z_{t_1}\ \cdots\ z_{t_c}] \bar{\gamma}_{\ell,t}$.  Combining the definitions of $\gamma_{\ell,t}$ and $\bar{\gamma}_{\ell,t}$ along with \eqref{main_assumption}, it follows for any such $t$ that
  \bequationNN
  \baligned
  &\gamma_{\ell,t}^T\bbmatrix 0_{\ell+1} & \cdots & 0_{\ell+1} \\ \psi_{\ell+1,t_1} & \cdots & \psi_{\ell+1,t_c} \ebmatrix^TL^TL\bbmatrix 0_{\ell+1} & \cdots & 0_{\ell+1} \\ \psi_{\ell+1,t_1} & \cdots & \psi_{\ell+1,t_c} \ebmatrix \\
   &\qquad \qquad= \ \bbmatrix 0_{\ell+1} \\ \psi_{\ell+1,t} \ebmatrix^TL^TL\bbmatrix 0_{\ell+1} & \cdots & 0_{\ell+1} \\ \psi_{\ell+1,t_1} & \cdots & \psi_{\ell+1,t_c} \ebmatrix \\
   &\qquad \qquad= \ z_t^T\bbmatrix z_{t_1} & \cdots & z_{t_c} \ebmatrix \\
   &\qquad \qquad= \ \bar{\gamma}_{\ell,t}^T\bbmatrix z_{t_1} & \cdots & z_{t_c} \ebmatrix^T\bbmatrix z_{t_1} & \cdots & z_{t_c} \ebmatrix \\
   &\qquad \qquad= \ \bar{\gamma}_{\ell,t}^T\bbmatrix 0_{\ell+1} & \cdots & 0_{\ell+1} \\ \psi_{\ell+1,t_1} & \cdots & \psi_{\ell+1,t_c} \ebmatrix^TL^TL\bbmatrix 0_{\ell+1} & \cdots & 0_{\ell+1} \\ \psi_{\ell+1,t_1} & \cdots & \psi_{\ell+1,t_c} \ebmatrix,
  \ealigned
  \eequationNN
  from which it follows that \rev{$\gamma_{\ell,t} = \bar{\gamma}_{\ell,t}$}.  Thus, with \eqref{eq.affine_1} and the definitions of \rev{$\gamma_{\ell,t}$} and \rev{$\bar{\gamma}_{\ell,t}$}, it follows for any $t\in\{\ell+1,\dots,m-1\}\backslash\{t_1,\dots,t_c\}$ that 
  \bequationNN
  \baligned
      \phi_t(a_{t,2})^T\phi_{\ell}(a_{\ell,2}^*) &=\ \bbmatrix 0_{\ell+1} \\ \psi_{\ell+1,t} \ebmatrix^TL^TL\bbmatrix 0_{\ell} \\ a_{\ell,2}^* + S_{\ell+1:m}^T(b_{\ell}y_0+y_{\ell}) \ebmatrix \\
    &=\ \gamma_{\ell,t}^T\bbmatrix 0_{\ell+1} & \cdots & 0_{\ell+1} \\ \psi_{\ell+1,t_1} & \cdots & \psi_{\ell+1,t_c} \ebmatrix^TL^TL\bbmatrix 0_{\ell} \\ a_{\ell,2}^* + S_{\ell+1:m}^T(b_{\ell}y_0+y_{\ell}) \ebmatrix \\
    &=\ \gamma_{\ell,t}^T\bbmatrix z_{t_1} & \cdots & z_{t_c} \ebmatrix^T z_{\ell} \\
    &=\ \bar{\gamma}_{\ell,t}^T\bbmatrix z_{t_1} & \cdots & z_{t_c} \ebmatrix^Tz_{\ell} = z_t^Tz_{\ell},
  \ealigned
  \eequationNN
  Combining this with \eqref{eq.affine_1}, it follows that $a_{\ell,2}^*$ from \eqref{beta} with $\beta_\ell$ from \eqref{beta_value} satisfies \eqref{eq.affine_begin}--\eqref{eq.affine_end}, as desired.  Let us show now that this $a_{\ell,2}^*$ also satisfies \eqref{starting_point}.  Indeed, by \eqref{beta}, \eqref{beta_value}, and \eqref{main_assumption}, it follows that
  \begin{align}
      &\phi_{\ell}(a_{\ell,2}^*)^T\phi_{\ell}(a_{\ell,2}^*) \nonumber \\ 
      &\ =\ \bbmatrix 0_{\ell} \\ a_{\ell,2}^* + S_{\ell+1:m}^T(b_{\ell}y_0 + y_{\ell}) \ebmatrix^TL^TL\bbmatrix 0_{\ell} \\ a_{\ell,2}^* + S_{\ell+1:m}^T(b_{\ell}y_0 + y_{\ell}) \ebmatrix \nonumber \\
    & \ =\ \beta_{\ell}^T\bbmatrix 0_{\ell+1} & \cdots & 0_{\ell+1} \\ \psi_{\ell+1,t_1} & \cdots & \psi_{\ell+1,t_c} \ebmatrix^TL^TL\bbmatrix 0_{\ell+1} & \cdots & 0_{\ell+1} \\ \psi_{\ell+1,t_1} & \cdots & \psi_{\ell+1,t_c} \ebmatrix \beta_{\ell} \nonumber \\
    & \ =\ z_{\ell}^T \bbmatrix z_{t_1}^T \\ \vdots \\ z_{t_c}^T \ebmatrix^T \(\bbmatrix 0_{\ell+1} & \cdots & 0_{\ell+1} \\ \psi_{\ell+1,t_1} & \cdots & \psi_{\ell+1,t_c} \ebmatrix^TL^TL\bbmatrix 0_{\ell+1} & \cdots & 0_{\ell+1} \\ \psi_{\ell+1,t_1} & \cdots & \psi_{\ell+1,t_c} \ebmatrix\)^{-1}\bbmatrix z_{t_1}^T \\ \vdots \\ z_{t_c}^T \ebmatrix z_{\ell} \nonumber \\
    & \ =\ z_{\ell}^T \bbmatrix z_{t_1}^T \\ \vdots \\ z_{t_c}^T \ebmatrix^T \(\bbmatrix z_{t_1}^T \\ \vdots \\ z_{t_c}^T \ebmatrix\bbmatrix z_{t_1} & \cdots & z_{t_c} \ebmatrix\)^{-1} \bbmatrix z_{t_1}^T \\ \vdots \\ z_{t_c}^T \ebmatrix z_{\ell} \leq z_{\ell}^Tz_{\ell}, \nonumber 
  \end{align}
  where the last inequality comes from the fact that the eigenvalues of
  \bequationNN
  \bbmatrix z_{\ell_1} & \cdots & z_{\ell_c} \ebmatrix \(\bbmatrix z_{\ell_1}^T \\ \vdots \\ z_{\ell_c}^T \ebmatrix\bbmatrix z_{\ell_1} & \cdots & z_{\ell_c} \ebmatrix\)^{-1} \bbmatrix z_{\ell_1}^T \\ \vdots \\ z_{\ell_c}^T \ebmatrix
  \eequationNN
  are all in $\{0,1\}$.  (As an aside, one finds that the inequality above is strict if $z_{\ell} \notin \linspan\{ z_{t_1} ,\dots,z_{t_c}\}$.)  Hence, we have shown that $a_{\ell,2}^*$ from \eqref{beta} satisfies \eqref{starting_point}.
  
  As previously mentioned (in the text following~\eqref{starting_point}), our goal now is to show that there exists a nonzero $\abar_{\ell,2}\in\R{m-\ell}$ such that $a_{\ell,2}^* + \lambda_{\ell}\abar_{\ell,2}$ satisfies \eqref{eq.affine_begin}--\eqref{eq.affine_end} for arbitrary $\lambda_{\ell}$.  From \eqref{eq.affine_begin}--\eqref{eq.affine_end}, such an $\abar_{\ell,2}\in\R{m-\ell}$ must satisfy
   \bequation\label{null.N}
   \bbmatrix 0_{\ell+1} & \cdots & 0_{\ell+1} \\ \psi_{\ell+1,\ell+1} & \cdots & \psi_{\ell+1,m-1} \ebmatrix^TL^TL\bbmatrix 0_{\ell} \\ \abar_{\ell,2} \ebmatrix = 0_{m-(\ell+1)}.
   \eequation
   Since 
   \bequation\label{eq.baoyu_matrix}
     \bbmatrix 0_{\ell+1} & \cdots & 0_{\ell+1} \\ \psi_{\ell+1,\ell+1} & \cdots & \psi_{\ell+1,m-1} \ebmatrix^TL^TL \in \R{(m-(\ell+1))\times m},
   \eequation
   it follows that this matrix has a null space of dimension at least $\ell+1$, i.e., there exist at least $\ell+1$ linearly independent vectors in $\R{m}$ belonging to the null space of this matrix.  Let $N_{\ell+1} \in \R{m \times (\ell+1)}$ be a matrix whose columns are $\ell+1$ linearly independent vectors in $\R{m}$ lying in the null space of \eqref{eq.baoyu_matrix}.  Since this null space matrix has $\ell+1$ columns, there exists a nonzero vector $\zeta_{\ell+1} \in \R{\ell+1}$ such that the first $\ell$ elements of $N_{\ell+1}\zeta_{\ell+1}$ are zero.  Letting 
  \bequationNN
  [\abar_{\ell,2}]_t := [N_{\ell+1}\zeta_{\ell+1}]_{\ell+t}\ \ \text{for all}\ \ t \in \{1,\dots,m-\ell\},
  \eequationNN   
one finds that $\abar_{\ell,2}$ satisfies \eqref{null.N}, as desired.  Consequently, as stated in the text following \eqref{starting_point}, it follows by the fact that the left-hand side of \eqref{eq.quadratic} is a strongly convex quadratic in the unknown~$\lambda_{\ell}$ and the fact that \eqref{starting_point} holds that we can claim that there exists $\lambda_{\ell}^* \in \R{}$ such that $a_{\ell,2} = a_{\ell,2}^* + \lambda_{\ell}^* \abar_{\ell,2} \in \R{m-\ell}$ satisfies \eqref{eq.linear_quadratic}.
  
  Combining all previous aspects of the proof, we have shown the existence of $A \in \R{m\times(m-1)}$ and $b\in\R{m-1}$ such that, with $\Ytilde_{1:m} \in \R{n \times m}$ defined as in \eqref{eq.Ytilde_formula}, the equations~\eqref{eq.key} hold.  The remaining desired conclusions, namely, that \eqref{eq.products_remain_positive} holds and that $\BFGS(W,S_{0:m},Y_{0:m}) = \BFGS(W,S_{1:m},\Ytilde_{1:m})$, follow from the existence of $\Ytilde_{1:m} \in \R{n \times m}$ (that we have proved), the fact that the equations in \eqref{eq.products_remain_positive} are a subset of the equations in \eqref{eq.R=Rtilde_simplified}, and the fact that \eqref{eq.key} were derived explicitly to ensure that, with $\Ytilde_{1:m}$ satisfying \eqref{eq.Ytilde_formula}, one would find that \eqref{eq.equal} holds. \qed
\end{proof}

\subsection{Implementing \AggBFGS}\label{eq.implementation}

We now discuss how one may implement our \AggBFGS{} scheme to iteratively aggregate displacement information in the context of an optimization algorithm employing BFGS approximations.  We also discuss the dominant computational costs of applying the scheme, and comment on certain numerical considerations that one should take into account in a software implementation.  The procedures presented in this section are guided by our proof of Theorem~\ref{th.real}.

As will become clear in our overall approach, in contrast to a traditional limited-memory scheme in which one always maintains the most recent curvature pairs to define a Hessian approximation, the pairs used in our approach might come from a subset of the previous iterations, with the gradient displacements potentially having been modified through our aggregation mechanism.  For concreteness, suppose that during the course of the run of an optimization algorithm for solving~\eqref{prob.f}, one has accumulated a set of curvature pairs, stored in the sets
\bequationNN
  \Scal := \{s_{k_0}, \dots, s_{k_{m-1}}\}\ \ \text{and}\ \ \Ycal := \{y_{k_0}, \dots, y_{k_{m-1}}\}
\eequationNN
where $\{k_i\}_{i=0}^{m-1} \subset \N{}$ with $k_i < k_{i+1}$ for all $i \in \{0,\dots,m-2\}$, such that the vectors in the former set (i.e., the iterate displacements) are linearly independent.  (Here, the elements of $\Ycal$ are not necessarily the gradient displacements computed in iterations $\{k_0,\dots,k_{m-1}\}$, but, for simplicity of notation, we denote them as $\Ycal$ even though they might have been modified during a previous application of our aggregation scheme.)  Then, suppose that a new curvature pair $(s_{k_m},y_{k_m})$ for $k_m \in \N{}$ with $k_{m-1} < k_m$ is available.  Our goal in this section is to show how one may add and, if needed, aggregate the information in these pairs in order to form new sets
\bequationNN
  \tilde\Scal \subseteq \Scal \cup \{s_{k_m}\}\ \ \text{and}\ \ \tilde\Ycal
\eequationNN
such that
\benumerate
  \item[(i)] the sets $\tilde\Scal$ and $\tilde\Ycal$ have the same cardinality, which is either $m$ or $m+1$,
  \item[(ii)] the vectors in the set $\tilde\Scal$ are linearly independent, and
  \item[(iii)] the BFGS inverse Hessian approximation generated from the data in $(\Scal \cup \{s_{k_m}\}, \Ycal \cup \{y_{k_m}\})$ is the same as the approximation generated from $(\tilde\Scal,\tilde\Ycal)$.
\eenumerate
As in the previous section, we henceforth simplify the subscript notation and refer to the ``previously stored'' displacement vectors as those in the sets $\{s_0,\dots,s_{m-1}\}$ and $\{y_0,\dots,y_{m-1}\}$, and refer to the ``newly computed'' curvature pair as $(s_m,y_m)$.

Once the newly computed pair $(s_m,y_m)$ is available, there are three cases.

\benumerate
  \item[Case 1.] The set of vectors $\{s_0,s_1,\dots,s_{m} \}$ is linearly independent. In this case, one simply adds the new curvature pair, which leads to the $(m+1)$-element sets
  \bequationNN
    \tilde\Scal = \{s_0,\dots,s_{m-1},s_m\}\ \ \text{and}\ \ \tilde\Ycal = \{y_0,\dots,y_{m-1},y_m\}.
  \eequationNN
  Notice that if $m=n$, then this case is not possible.
  \item[Case 2.] The new iterate displacement vector $s_m$ is parallel to the most recently stored iterate displacement vector, i.e., $s_{m-1} = \tau s_m$ for some $\tau \in \R{}$.  In this case, one should discard the most recently stored pair and replace it with the newly computed one, which leads to the $m$-element sets
  \bequationNN
    \tilde\Scal = \{s_0,\dots,s_{m-2},s_m\}\ \ \text{and}\ \ \tilde\Ycal = \{y_0,\dots,y_{m-2},y_m\}.
  \eequationNN
  This choice is justified by Theorem~\ref{th.1_update_skip}.
  \item[Case 3.] For some $j \in \{0,\dots,m-2\}$, an iterate displacement vector $s_j$ lies in the span of the subsequent iterate displacements vectors, i.e., $s_j \in \linspan \{s_{j+1},\dots,s_m\}$.  In this case, one should apply our aggregation scheme to determine the vectors $\{\ytilde_{j+1},\dots,\ytilde_m\}$, then remove the pair $(s_j,y_j)$, leading to the $m$-element sets
  \bequationNN
    \tilde\Scal = \{s_0,\dots,s_{j-1},s_{j+1},\dots,s_m\}\ \ \text{and}\ \ \tilde\Ycal = \{y_0,\dots,y_{j-1},\ytilde_{j+1},\dots,\ytilde_m\}.
  \eequationNN
  This choice is justified by Theorem~\ref{th.real}.
\eenumerate

Computationally, the first step is to determine which of the three cases occurs.  One way to do this efficiently is to maintain a Cholesky factorization of an inner product matrix corresponding to the previously stored iterate displacement vectors, then attempt to add to it a new row/column corresponding to the newly computed iterate displacement, checking whether the procedure breaks down.  Specifically, before considering the newly computed pair $(s_m,y_m)$, suppose that one has a lower triangular matrix $\Theta \in \R{m \times m}$ with positive diagonal elements such that
\bequation\label{eq.fact_m-1}
  \bbmatrix s_{m-1} & \cdots & s_0 \ebmatrix^T\bbmatrix s_{m-1} & \cdots & s_0 \ebmatrix = \Theta\Theta^T,
\eequation
which exists due to the fact that the vectors in $\{s_0,\dots,s_{m-1}\}$ are linearly independent.  A Cholesky factorization of an augmented inner product matrix would consist of a scalar $\mu \in \R{}_{>0}$, vector $\delta \in \R{m}$, and lower triangular $\Delta \in \R{m \times m}$ with
\bequation\label{eq.fact_m}
  \bbmatrix s_m & s_{m-1} & \cdots & s_0 \ebmatrix^T\bbmatrix s_m & s_{m-1} & \cdots & s_0 \ebmatrix = \bbmatrix \mu & 0 \\ \delta & \Delta \ebmatrix \bbmatrix \mu & \delta^T \\ 0 & \Delta^T \ebmatrix.
\eequation
As is well known, equating terms in \eqref{eq.fact_m-1} and \eqref{eq.fact_m} one must have $\mu = \|s_m\|$, $\delta^T = [s_m^Ts_{m-1}\ \cdots\ s_m^Ts_0]/\mu$, and $\Delta\Delta^T = \Theta\Theta^T - \delta\delta^T$, meaning that $\Delta$ can be obtained from~$\Theta$ through a \emph{rank-one downdate}; e.g., see \cite{GillGoluMurrSaun74}.  If this downdate does not break down---meaning that all diagonal elements of $\Delta$ are computed to be positive---then one is in Case~1 and the newly updated Cholesky factorization has been made available when yet another curvature pair is considered (after a subsequent optimization algorithm iteration).  Otherwise, if the downdate does break down, then it is due to a computed diagonal element being equal to zero.  This means that, for some smallest $i \in \{1,\dots,m\}$, one has found a lower triangular matrix $\Xi \in \R{i \times i}$ with positive diagonal elements and a vector $\xi \in \R{i}$ such that
\bequation\label{eq.downdate_fail}
  \bbmatrix s_m & s_{m-1} & \cdots & s_{m-i} \ebmatrix^T\bbmatrix s_m & s_{m-1} & \cdots & s_{m-i} \ebmatrix = \bbmatrix \Xi & 0 \\ \xi^T & 0 \ebmatrix \bbmatrix \Xi^T & \xi \\ 0 & 0 \ebmatrix.
\eequation
Letting $\tau \in \R{i}$ be the unique vector satisfying $\Xi^T\tau = \xi$, one finds that the vector \rev{$[\tau^T,-1]^T$} lies in the null space of \eqref{eq.downdate_fail}, from which it follows that
\bequationNN
  \bbmatrix s_m & s_{m-1} & \cdots & s_{m-i+1} \ebmatrix \tau = s_{m-i},
\eequationNN
where the first element of $\tau$ must be nonzero since $\{s_{m-1},\dots,s_{m-i}\}$ is a set of linearly independent vectors.  If the breakdown occurs for $i=1$, then one is in Case~2.  If the breakdown occurs for $i > 1$, then one is in Case~3 with the vector $\tau$ that one needs to apply our aggregation scheme to remove the pair $(s_{m-i},y_{m-i})$.

Notice that if the breakdown occurs in the rank-one downdate as described in the previous paragraph, then one can continue with standard Cholesky factorization updating procedures in order to have the factorization of
\bequationNN
  \bbmatrix s_m & \cdots & s_{m-i+1} & s_{m-i-1} & \cdots s_0 \ebmatrix^T\bbmatrix s_m & \cdots & s_{m-i+1} & s_{m-i-1} & \cdots s_0 \ebmatrix
\eequationNN
available in subsequent iterations.  For brevity and since it is outside of our main scope, we do not discuss this in detail.  Overall, the computational costs so far are $\Ocal(mn)$ (for computing the inner products $\{s_m^Ts_{m-1},\dots,s_m^Ts_0\}$) plus $\Ocal(m^2)$ (for updating the Cholesky factorization and, in Case~2 or Case~3, computing $\tau$).

If Case~1 occurs, then no additional computation is necessary; one simply adds the new curvature pair as previously described.  Similarly, if Case~2 occurs, then again no additional computation is necessary; one simply replaces the most recent stored pair with the newly computed one.  Therefore, we may assume for the remainder of this section that Case~3 occurs, we have identified an index $j$ ($=m-i$ using the notation above) such that $s_j \in \linspan\{s_{j+1},\dots,s_m\}$, and we have $\tau \in \R{m-j}$ such that $s_j = S_{j+1:m}\tau$.  Our goal then is to apply our aggregation scheme to modify the gradient displacement vectors $\{y_{j+1},\dots,y_m\}$ to compute
\bequation\label{eq.newYtilde}
  \Ytilde_{j+1:m} = W_{0:j-1}^{-1}S_{j+1:m} \bbmatrix A & 0 \ebmatrix + y_j \bbmatrix b \\ 0 \ebmatrix^T + Y_{j+1:m},
\eequation
where $W_{0:j-1}$ represents the BFGS inverse Hessian approximation defined by some initial positive definite matrix $W \succ 0$ and the curvature pairs $\{s_i,y_i\}_{i=0}^{j-1}$, and where $A \in \R{(m-j) \times (m-j-1)}$ and $b\in \R{m-j-1}$ are the unknowns to be determined.

For simplicity in the remainder of this section, let us suppose that $j=0$ so that the pair $(s_0,y_0)$ is to be removed as in the notation of~\S\ref{eq.basics_aggBFGS} and \S\ref{sec.existence}.  As one might expect, this is the value of $j$ for which the computational costs of computing $A \in \R{m \times (m-1)}$ and $b \in \R{m-1}$ are the highest.  For all other values of~$j$, there is a cost for computing $W_{0:j-1}^{-1}S_{j+1:m}$ as needed in \eqref{eq.newYtilde}.  This matrix can be computed \emph{without} forming the BFGS Hessian approximation $W_{0:j-1}^{-1}$; it can be constructed, say, by computing matrix-vector products with a compact representation of this approximation for a total cost of $\rev{\Ocal(j(m-j)n)} \leq \Ocal(m^2n)$; see \S7.2 in~\cite{NoceWrig06}.

Let us now describe how one may implement our aggregation scheme such that, given $S_{1:m}$ with full column rank, $Y_{1:m}$, $\rho_{1:m} > 0$, $\tau \in \R{m}$ satisfying $s_0 = S_{1:m}\tau$, $y_0$, and $\rho_0 > 0$, one may compute $A \in \R{m \times (m-1)}$ and $b \in \R{m-1}$ in order to obtain $\Ytilde_{1:m}$ as in \eqref{eq.Ytilde_formula}.  By Theorem~\ref{th.real}, it follows that real values for $A$ and~$b$ exist to satisfy~\eqref{eq.key}.  The computation of the vector $b$ is straightforward; it can be computed by the formula \eqref{eq.b_simplified}.  Assuming that the products in $S_{1:m}^TY_{1:m-1}$ have already been computed (\rev{which, using previously computed inner products and recomputing any as needed if/when gradient displacements have been modified by aggregation, costs $\Ocal((m-j)^2n)$}), the cost of computing $b$ is $\Ocal(m^2)$.

For computing~$A$, let us first establish some notation since the elements of this matrix will be computed with a specific order.  As in the proof of Theorem~\ref{th.real}, let $A = [a_1\ \cdots\ a_{m-1}]$ where $a_\ell \in \R{m}$ for all $\ell \in \{1,\dots,m-1\}$, and, as in \eqref{eq.aj}, let
\bequationNN
  a_\ell = Q^{-1}\bbmatrix a_{\ell,1} \\ a_{\ell,2} \ebmatrix,\ \text{where}\ a_{\ell,1} \in \R{\ell},\ a_{\ell,2} \in \R{m-\ell},\ \text{and}\ Q := S_{1:m}^TW^{-1}S_{1:m} \succ 0.
\eequationNN
Here and going forward, our computations require products with $Q^{-1}$.  Rather than form this matrix explicitly, one may maintain a Cholesky factorization of~$Q$ and add/delete rows/columns---as described previously for the inner product matrix corresponding to the iterate displacements---as the iterate displacement set is updated throughout the run of the (outer) optimization algorithm.  With such a factorization, products with $Q^{-1}$ are obtained by triangular solves in a standard fashion.  If $W \succ 0$ is diagonal or defined by a limited memory approximation, then the cost of updating this factorization in each instance is $\Ocal(mn)$ (for computing $W^{-1}s_m$) plus $\Ocal(m^2)$ for updating the factorization.  Each backsolve costs $\Ocal(m^2)$.

\rev{An approach for computing~$A$ is now summarized as Algorithm~\ref{alg.A}.  As explained in the proof of Theorem~\ref{th.real}, the values of $\{a_{\ell,1}\}_{\ell=1}^{m-1}$ are set to ensure that \eqref{eq.R=Rtilde_simplified} is satisfied.  Assuming that $S_{1:m}^Ty_0$ has been computed at cost $\Ocal(mn)$, the cost of computing $\{a_{\ell,1}\}_{\ell=1}^{m-1}$ is $\Ocal(m^2)$.  The values for $\{a_{\ell,2}\}_{\ell=1}^{m-1}$ are then computed in reverse order to solve the system of linear and quadratic equations comprising \eqref{eq.YWYtilde_simplified}.  Assuming (as has been mentioned) that a factorization of $Q$ is available, and assuming that the elements of the right-hand side of \eqref{eq.YWYtilde_compact_reduced} have been pre-computed (at cost $\Ocal(m^3)$), the cost of computing $a_{m-1,2}$ is $\Ocal(1)$.  As for computing the remaining vectors, the presented scheme follows the proof of Theorem~\ref{th.real}.  The most expensive operation in each iteration of this scheme is the QR factorization of the matrix in \eqref{eq.baoyu_matrix}, which for each $\ell$ is of size $(m - (\ell+1)) \times m$.  Summing the cost of these for $\ell=m-2$ to $\ell=1$, the total cost is found to be $\Ocal(m^4)$.}

\balgorithm[ht]
\rev{
  \caption{:\ \rev{Displacement Aggregation, Computation of $A$}}
  \label{alg.A}
  \balgorithmic[1]
    \State For each $\ell \in \{1,\dots,m-1\}$, compute the $\ell$-element vector $a_{\ell,1}$ by~\eqref{eq.aj1}.
    \State Compute $a_{m-1,2}$ by solving the quadratic equation \eqref{eq.quadratic}.
    \For{$\ell = m-2,\dots,1$}
      \State Compute $a^*_{\ell,2}$ to satisfy \eqref{beta}--\eqref{beta_value}.
      \State Compute nonzero vector $\bar{a}_{\ell,2}$ to satisfy \eqref{null.N}.
      \State Compute $\lambda_\ell^* \in \R{}$ such that $a_{\ell,2} = a^*_{\ell,2} + \lambda_\ell^* \bar{a}_{\ell,2}$ solves the quadratic equation~\eqref{eq.quadratic}.
    \EndFor
    \State \Return $A = \bbmatrix a_1 & \cdots & a_{m-1} \ebmatrix$ with $\{a_j\}_{j=1}^{m-1}$ defined by \eqref{eq.aj}.
  \ealgorithmic
}
\ealgorithm

This completes our description of a manner in which our \AggBFGS{} scheme can be implemented.  Observing the computational costs that have been cited, one finds that \rev{a conservative estimate of the total cost is $\Ocal(m^2n) + \Ocal(m^4)$.  Upon closer inspection, one finds that for $j \approx 0$ the cost is dominated by the cost of computing/updating $S_{1:m}^T\Ytilde_{1:m}$ \emph{after} aggregation has been performed, whereas for $j \approx m$ the cost is dominated by the cost of computing $W_{0:j-1}^{-1}S_{j+1:m}$ \emph{before} aggregation is performed.  In either case, if $m \ll n$ (say with $3 \leq m \leq 10$, as is typical in practice), then the overall cost is not too dissimilar from $\Ocal(4mn)$ per iteration, which is the cost of other limited-memory schemes such as L-BFGS \cite{ByrdNoceSchn94}.}

We end this section by stating the following result, for closure.

\btheorem
  If one applies \AggBFGS{} as described in this section, then one need only store at most $m \leq n$ curvature pairs such that the corresponding BFGS inverse Hessian approximations are equivalent to those in a full-memory BFGS scheme.  Consequently, under the same conditions as in  Theorem~\ref{th.bfgs_superlinear}, the resulting optimization algorithm produces $\{x_k\}$ that converges to $x_*$ at a superlinear rate.
\etheorem

\section{Numerical Demonstrations}\label{sec.numerical}

Our goals in this section are to provide additional numerical demonstrations of applications of \AggBFGS{}.  (Recall that a preview demonstration has been provided in Figure~\ref{fig.lbfgs_vs_aggbfgs}.)  Our first goal is to show empirically that limited-memory-type BFGS inverse Hessian approximations provided by \AggBFGS{} accurately represent the approximations provided by full-memory BFGS.  We demonstrate that while numerical errors might accumulate to some degree as \AggBFGS{} is applied over a sequence of iterations, the inverse Hessian approximations provided by \AggBFGS{} are not necessarily poor after multiple iterations.  Our second goal is to provide the results of numerical experiments with an adaptive L-BFGS method that uses our \AggBFGS{} scheme and show that over a diverse test set it can outperform a standard L-BFGS approach.  These experiments are run in the practical regime when the number of pairs stored and employed is small relative to $n$.

\subsection{Simulated Data}

We implemented \AggBFGS{} in MATLAB and ran two sets of experiments.  First, for $(n,m) \in \{4,8,16,32,64,128\}^2$ with $m \leq n$, we generated data to show the error of applying \AggBFGS{} to aggregate a single curvature pair.  For each $(n,m)$, we generated 100 datasets using the following randomized procedure.  First, MATLAB's built-in \texttt{sprandsym} function was used to generate a random positive definite matrix with condition number approximately $10^4$.  This matrix defined a quadratic function.  Second, a random fixed point was determined using MATLAB's built-in \texttt{randn} routine.  From this point, a mock optimization procedure for minimizing the generated quadratic was run to generate $\{(s_k,y_k)\}_{k=1}^m$; in particular, for $k \in \{1,\dots,m\}$, starting with the fixed point, a descent direction was chosen as the negative gradient plus noise (specifically, the norm of the gradient divided by 10 times a random vector generated with \texttt{randn}) and a stepsize was chosen by an exact line search to compute the subsequent iterate and gradient displacement pair.  Third, a vector $\tau \in \R{m}$ was generated randomly using \texttt{randn} in order to define $s_0 = [s_1\ \cdots\ s_m]\tau$.  The corresponding gradient displacement~$y_0$ was determined by stepping \emph{backward} from the fixed point along~$s_0$.  In this manner, we obtained $\{(s_k,y_k)\}_{k=0}^m$ from a mock optimization procedure from some initial point in such a way that $s_0$ was guaranteed to lie the span of the subsequent iterate displacements.

\begin{figure}[ht]
  \centering
  \includegraphics[width=\textwidth,clip=true,trim=170 0 150 20]{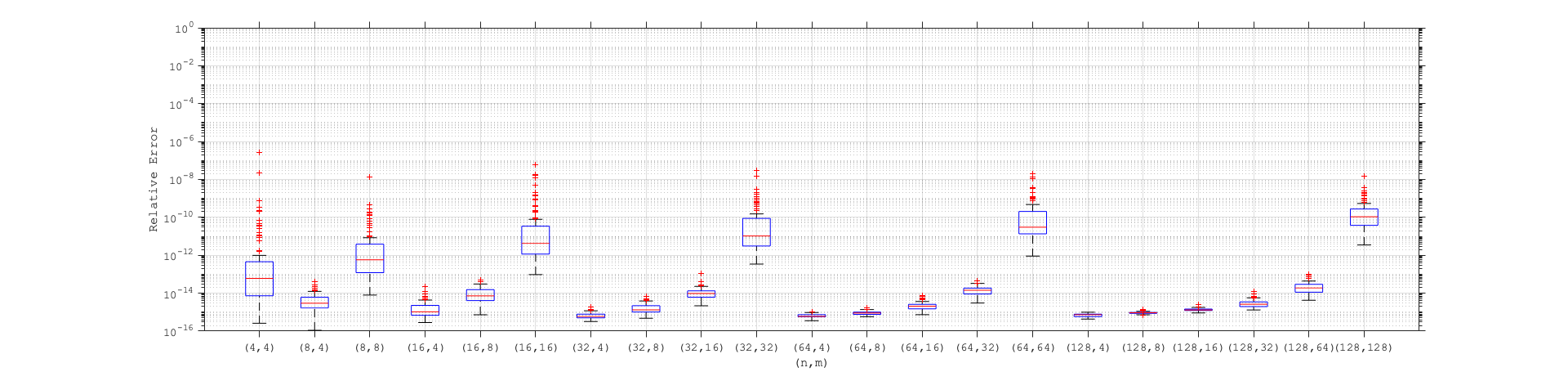}
  \caption{Relative errors of the differences between BFGS inverse Hessian approximations (computed from $\{(s_k,y_k)\}_{k=0}^m$) and the corresponding inverse Hessian approximations represented after applying \AggBFGS{} to aggregate the information from a single curvature pair (into $\{(s_k,\ytilde_k)\}_{k=1}^m$).  For each $(n,m)$, the box plots show the results for 100 randomly generated instances.  Relative error is the maximum absolute difference between corresponding matrix entries divided by the largest absolute value of an element of the BFGS inverse Hessian approximation.}
  \label{fig.one_aggregation_errors}
\end{figure}

For each dataset starting from $W_0 = I$, we computed the BFGS inverse Hessian approximation from $\{(s_k,y_k)\}_{k=0}^m$ and constructed the inverse Hessian approximation corresponding to the set of pairs when \AggBFGS{} was used to aggregate the information into $\{(s_k,\ytilde_k)\}_{k=1}^m$.  Figure~\ref{fig.one_aggregation_errors} shows box plots for relative errors between each pair of inverse Hessian approximations, where error is measured in terms of the maximum component-wise absolute difference between matrix entries divided by the largest absolute value of an element of the BFGS inverse Hessian approximation.  The results show that while the errors are larger as $n$ increases and as $m$ is closer to $n$, they remain accurate relative to machine precision for all $(n,m)$.

\begin{figure}[ht]
  \centering
  \begin{subfigure}[b]{0.31\textwidth}
    \includegraphics[width=\textwidth,clip=true,trim=0 0 40 15]{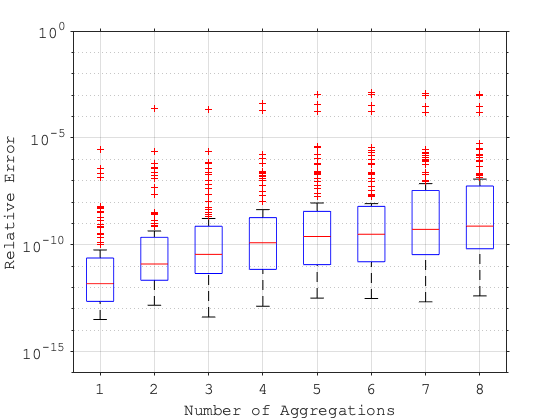}
    \caption{$n=m=8$}
  \end{subfigure}
  ~
  \begin{subfigure}[b]{0.31\textwidth}
    \includegraphics[width=\textwidth,clip=true,trim=0 0 40 15]{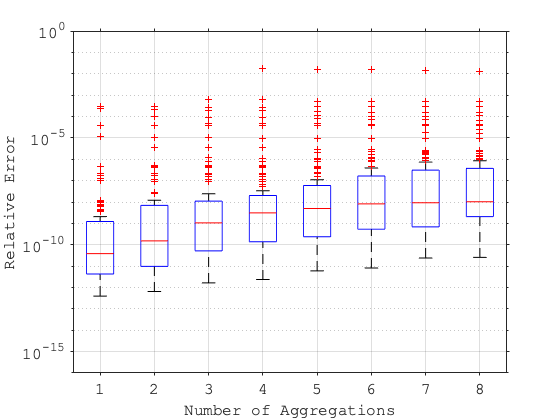}
    \caption{$n=m=32$}
  \end{subfigure}
  ~
  \begin{subfigure}[b]{0.31\textwidth}
    \includegraphics[width=\textwidth,clip=true,trim=0 0 40 15]{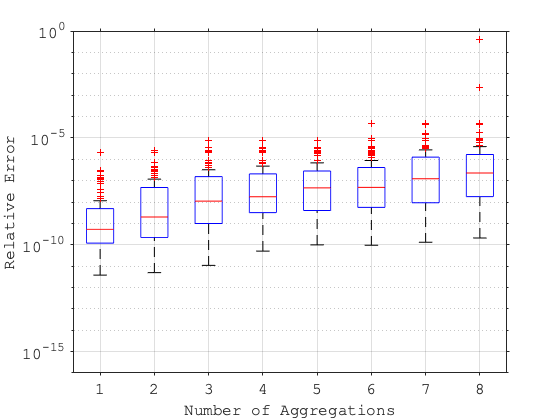}
    \caption{$n=m=128$}
  \end{subfigure}
  \caption{Accumulation of relative errors of the differences between BFGS inverse Hessian approximations and the corresponding inverse Hessian approximations after applying \AggBFGS{} to aggregate the information.  For each $(n,m)$, the box plots show the results for 100 randomly generated instances.  Relative error is the maximum absolute difference between corresponding matrix entries divided by the largest absolute value of an element of the BFGS inverse Hessian approximation.}
  \label{fig.multiple_aggregation_errors}
\end{figure}

As a second experiment, for $(n,m) \in \{(8,8),(32,32),(128,128)\}$, we aimed to investigate how errors might accumulate as \AggBFGS{} is applied over a sequence of iterations.  For these experiments, we generated data using a similar mock optimization procedure as in the previous experiment.  Specifically, from some randomly generated starting point, we computed a random step as in the aforementioned procedure, continuing until $n+8$ iterations were performed.  Figure~\ref{fig.multiple_aggregation_errors} shows the errors for each iteration beyond $k=n$.  As in Figure~\ref{fig.lbfgs_vs_aggbfgs}, these results show that the errors do not accumulate too poorly as $k$ increases.  We conjecture that part of the reason for this is that errors that result from each application of \AggBFGS{} can be over-written eventually, at least to some extent, similar to the manner in which curvature information is ultimately over-written in full-memory BFGS.

\subsection{A practical adaptive L-BFGS method using \AggBFGS{}}

To demonstrate the use of \AggBFGS{} in a minimization algorithm, \rev{we} compare the results of \rev{three} algorithms using quasi-Newton inverse Hessian approximations for computing the search directions, where for each algorithm the same weak Wolfe line search was used for computing stepsizes.  The first algorithm employed L-BFGS approximations.  The second algorithm also employed L-BFGS approximations, but with \AggBFGS{} used to aggregate information when deemed appropriate in the manner described \rev{below.  Since this second algorithm does not always aggregate the oldest pair, but rather uses a particular scheme for choosing which historical pair to use in the aggregation approach (as explained later on), we also considered a third algorithm.  This algorithm uses the same scheme as the second to determine which pair to consider, but simply removes this pair rather than perform aggregation with it.  Consequently, the comparison between the second and third algorithms shows the effect of aggregation itself, rather than only the difference between removing the oldest versus some other historical pair in an L-BFGS scheme.}

\rev{We performed experiments with problems from the CUTEst collection \cite{GoulOrbaToin15}.}  We chose all problems from the collection for which \rev{the number of variables could be chosen} in the interval $[10,3000]$.  For \rev{most} problems, $n \gg m$, meaning that the costs of performing aggregation was negligible compared to the costs of computing search directions, which were the same for \rev{all} algorithms using the standard two-loop recursion for L-BFGS.  \rev{(See Tables~\ref{tab.results1} and \ref{tab.results2} for $n$ for all test problems used.)}  \rev{All} algorithms were run until an iteration, call it $k \in \N{}$, was reached with $\|g_k\|_\infty \leq 10^{-6} \max\{1,\|g_0\|_\infty\}$, or until an iteration limit of $10^5$ was surpassed.

For the second algorithm, \AggBFGS{} was employed in the following manner.  Suppose in iteration $k \in \N{}$ that the algorithm had in storage the displacement pairs $\{(s_{k_j},\ytilde_{k_j})\}_{j=1}^{\mbar}$ for some indices $\{k_1,\dots,k_{\mbar}\} \subseteq \{1,\dots,k-1\}$ with $\mbar \leq m$ and the new pair $(s_k,y_k)$ has been computed.  Using techniques described in~\S\ref{eq.implementation}, starting with $j=\mbar$ and iterating down to $j=1$, we determined if $s_{k_j}$ lay approximately in $\linspan\{s_{k_{j+1}},\dots,s_{k_{\mbar}},s_k\}$ (see next paragraph).  If so, then we applied \AggBFGS{} to aggregate the information in the pair with index $k_j$; otherwise, if no such iterate displacement was determined, then either the pair $(s_k,y_k)$ was simply added to the set of pairs in storage (if $\mbar < m$) or the pair with index $k_1$ was dropped (if $\mbar = m$, as in standard L-BFGS).  In this manner, the number of pairs stored and employed remained less than or equal to $m$, as usual for L-BFGS.

To determine if $s_{k_j}$ lay approximately in $\linspan\{s_{k_{j+1}},\dots,s_{k_{\mbar}},s_k\}$, we first computed the orthogonal projection of $s_{k_j}$ onto this subspace, call it $\shat_{k_j}$.  Computing this projection is inexpensive given a Cholesky factorization of the inner product matrix $[s_{k_{j+1}}\ \cdots\ s_{k_{\mbar}}\ s_k]^T[s_{k_{j+1}}\ \cdots\ s_{k_{\mbar}}\ s_k]$, which can be updated with only $\Ocal(mn)$ cost in each outer iteration as described in \S\ref{eq.implementation}.  If the condition
\bequation\label{eq.proj_error}
  \|s_{k_j} - \shat_{k_j}\|_2 \leq \texttt{tol} \cdot \|\shat_{k_j}\|_2
\eequation
held for some $\texttt{tol} \in (0,1)$, then it was determined that $s_{k_j}$ lay approximately in the subspace, since in this case the norm of the projection $\shat_{k_j}$ was sufficiently large compared to the orthogonal component $s_{k_j} - \shat_{k_j}$.  For $j = \{\mbar,\dots,2\}$, we used $\texttt{tol} = 10^{-8}$, while for $j = 1$ we loosened the tolerance to \rev{$\texttt{tol} = 10^{-4}$} to promote aggregation of the oldest pair.  If \eqref{eq.proj_error} held for some index $k_j$, then aggregation was performed with the pair $(\shat_{k_j},\ytilde_{k_j})$.  \rev{Notice that this ensures that the aggregation is performed with an iterate displacement vector (namely, the projected vector~$\shat_{k_j}$) that lies in the subspace defined by subsequent displacements, as required in our algorithms and theoretical results in \S\ref{sec.aggregation}.}  \rev{(We remark in passing that, using this scheme, it is possible that aggregation could be triggered multiple times in a single iteration.  Our software allows for this, but it did not occur in our experiments.)}

One might expect that, under these conditions, \AggBFGS{} might not perform aggregation often.  However, on the contrary, our results show that aggregation was performed quite often.\footnote{This provides evidence for the belief, held by some optimization researchers, that when solving certain large-scale problems one often observes that consecutive steps lie approximately in low-dimensional subspaces.}  \rev{Tables~\ref{tab.results1} and \ref{tab.results2} show} detailed results for \rev{the first two algorithms in our experiments for $m=5$.  (We do not include the third algorithm since, as shown below, it was inferior to the first and second algorithms.)}  For conciseness, we refer to the algorithm that employed a standard LBFGS(5) strategy as \texttt{LBFGS(5)}, and we refer to the strategy that employed \AggBFGS{} as \texttt{AggBFGS(5)}.  For \texttt{AggBFGS(5)}, the table reports the number of aggregations performed, which in many cases was significant compared to the number of iterations.

\btable[ht]
  \caption{Numbers of iterations, function evaluations, and aggregations when algorithms are applied to solve problems from the CUTEst set with $n \in [10,3000]$.}
  \label{tab.results1}
  \centering
  \tiny
  \texttt{
  \btabular{|lr|rrr|rr|}
    \hline
    & & \multicolumn{3}{c|}{AggBFGS(5)} & \multicolumn{2}{c|}{LBFGS(5)} \\
    Name & $n$ & Iters. & Funcs. & Aggs. & Iters. & Funcs. \\
    \hline
            ARGLINA & 200 &               3 &              12 &               1 &               3 &              12\\
        ARGLINB & 200 &               3 &               3 &               1 &               3 &               3\\
        ARGLINC & 200 &               3 &               3 &               1 &               3 &               3\\
      ARGTRIGLS & 200 &            1494 &           10220 &               0 &            1494 &           10220\\
        ARWHEAD & 1000 &               2 &               2 &               0 &               2 &               2\\
        BA-L1LS & 57 &              60 &             218 &               0 &              60 &             218\\
      BA-L1SPLS & 57 &              62 &             248 &               0 &              62 &             248\\
        BDQRTIC & 1000 &             160 &             734 &              28 &             399 &            2051\\
            BOX & 1000 &              16 &              76 &               6 &             --- &             ---\\
       BOXPOWER & 1000 &              20 &              86 &              11 &              12 &              29\\
        BROWNAL & 200 &               3 &               4 &               0 &               3 &               4\\
     BROYDN3DLS & 1000 &              39 &              45 &               0 &              39 &              45\\
       BROYDN7D & 1000 &            1444 &            2996 &               0 &            1444 &            2996\\
     BROYDNBDLS & 1000 &              71 &             151 &               0 &              71 &             151\\
         BRYBND & 1000 &              71 &             151 &               0 &              71 &             151\\
       CHAINWOO & 1000 &             563 &            2791 &              11 &             506 &            2577\\
       CHNROSNB & 50 &             186 &             355 &               0 &             186 &             355\\
       CHNRSNBM & 50 &             533 &            1593 &               0 &             533 &            1593\\
         COSINE & 1000 &             112 &             467 &               0 &             112 &             467\\
       CRAGGLVY & 1000 &             247 &            1189 &              10 &             241 &            1153\\
        CURLY10 & 1000 &             --- &             --- &             --- &             --- &             ---\\
        CURLY20 & 1000 &             --- &             --- &             --- &             --- &             ---\\
        CURLY30 & 1000 &             --- &             --- &             --- &             --- &             ---\\
       DIXMAANA & 300 &              16 &              48 &              10 &              20 &              60\\
       DIXMAANB & 300 &              22 &              78 &               6 &              24 &              97\\
       DIXMAANC & 300 &              20 &             106 &               5 &              18 &              78\\
       DIXMAAND & 300 &              25 &              93 &               6 &              35 &             135\\
       DIXMAANE & 300 &             400 &            1865 &               5 &             378 &            1798\\
       DIXMAANF & 300 &             383 &            1906 &               6 &             332 &            1728\\
       DIXMAANG & 300 &             264 &            1428 &               8 &             319 &            1526\\
       DIXMAANH & 300 &             425 &            2148 &              34 &             526 &            2668\\
       DIXMAANI & 300 &            1763 &            9092 &              33 &            2054 &           10878\\
       DIXMAANJ & 300 &             307 &            1380 &               6 &             406 &            2006\\
       DIXMAANK & 300 &             333 &            1589 &              15 &             331 &            1618\\
       DIXMAANL & 300 &             302 &            1553 &              29 &             363 &            1860\\
       DIXMAANM & 300 &            2901 &           14233 &              22 &            5220 &            6061\\
       DIXMAANN & 300 &             701 &            3749 &              15 &            1623 &            8508\\
       DIXMAANO & 300 &             694 &            3423 &              48 &             720 &            3642\\
       DIXMAANP & 300 &             555 &            2662 &              52 &             677 &            3421\\
       DIXON3DQ & 1000 &            5247 &            5978 &               0 &            5247 &            5978\\
     DMN15102LS & 66 &             206 &             962 &              44 &            1835 &           13002\\
     DMN15103LS & 99 &             853 &            4535 &             111 &            1696 &            9250\\
     DMN15332LS & 66 &            1253 &            6574 &             194 &             978 &            5241\\
     DMN15333LS & 99 &           11244 &           39107 &              26 &            4339 &           16755\\
     DMN37142LS & 66 &             639 &            3314 &             105 &             710 &            3744\\
     DMN37143LS & 99 &              12 &             128 &               0 &              12 &             128\\
        DQDRTIC & 1000 &              40 &             191 &              21 &              23 &              75\\
         DQRTIC & 1000 &             311 &            1838 &               6 &             511 &            2973\\
        EDENSCH & 36 &             110 &             547 &              11 &             138 &             673\\
            EG2 & 1000 &               5 &               5 &               3 &               5 &               5\\
       EIGENALS & 110 &             661 &            1898 &               0 &             661 &            1898\\
       EIGENBLS & 110 &            1239 &            3765 &               0 &            1239 &            3765\\
       EIGENCLS & 462 &            1898 &            5431 &               0 &            1898 &            5431\\
        ENGVAL1 & 1000 &              50 &              74 &               1 &              50 &              74\\
       ERRINROS & 50 &             900 &            4597 &              92 &            3648 &           16633\\
       ERRINRSM & 50 &             662 &            3225 &              60 &             --- &             ---\\
       EXTROSNB & 1000 &             301 &            1231 &               1 &             307 &            1350\\
       FLETBV3M & 1000 &              92 &             451 &               3 &             114 &             518\\
       FLETCBV2 & 1000 &            1014 &            1022 &               0 &            1014 &            1022\\
       FLETCBV3 & 1000 &             --- &             --- &             --- &             --- &             ---\\
       FLETCHBV & 1000 &             --- &             --- &             --- &             --- &             ---\\
       FLETCHCR & 1000 &             176 &            1625 &               0 &             176 &            1625\\
       FMINSRF2 & 961 &             374 &             429 &               0 &             374 &             429\\
       FMINSURF & 961 &             245 &             309 &               0 &             245 &             309\\
       FREUROTH & 1000 &              23 &              64 &               0 &              23 &              64\\
       GENHUMPS & 1000 &            3368 &            9766 &               0 &            3368 &            9766\\
        GENROSE & 500 &            2299 &           13271 &               0 &            2299 &           13271\\
       HILBERTA & 10 &              31 &             145 &               9 &              55 &             240\\
    \hline
  \etabular
  }
\etable

\btable[ht]
  \caption{Numbers of iterations, function evaluations, and aggregations when algorithms are applied to solve problems from the CUTEst set with $n \in [10,3000]$.}
  \label{tab.results2}
  \centering
  \tiny
  \texttt{
  \btabular{|lr|rrr|rr|}
    \hline
    & & \multicolumn{3}{c|}{AggBFGS(5)} & \multicolumn{2}{c|}{LBFGS(5)} \\
    Name & $n$ & Iters. & Funcs. & Aggs. & Iters. & Funcs. \\
    \hline
           HILBERTB & 50 &              14 &              38 &               6 &              14 &              38\\
       HYDC20LS & 99 &            1802 &            4128 &               0 &            1802 &            4128\\
          INDEF & 1000 &             --- &             --- &             --- &             --- &             ---\\
         INDEFM & 1000 &             --- &             --- &             --- &             --- &             ---\\
      INTEQNELS & 102 &              12 &              13 &               0 &              12 &              13\\
         JIMACK & 81 &            1584 &            3833 &               0 &            1584 &            3833\\
        LIARWHD & 1000 &              21 &              96 &              18 &              19 &              78\\
     LUKSAN11LS & 100 &            1332 &            5926 &               0 &            1332 &            5926\\
     LUKSAN12LS & 98 &             290 &             639 &               0 &             290 &             639\\
     LUKSAN13LS & 98 &              73 &             147 &               0 &              73 &             147\\
     LUKSAN14LS & 98 &             224 &            1041 &               0 &             224 &            1041\\
     LUKSAN15LS & 100 &              40 &             116 &               0 &              40 &             116\\
     LUKSAN16LS & 100 &              44 &              88 &               0 &              44 &              88\\
     LUKSAN17LS & 100 &             402 &            1917 &               0 &             402 &            1917\\
     LUKSAN21LS & 100 &             702 &            1964 &               0 &             702 &            1964\\
     LUKSAN22LS & 100 &             224 &             852 &               0 &             224 &             852\\
        MANCINO & 100 &              44 &             182 &              11 &              58 &             257\\
      MNISTS0LS & 494 &               2 &               2 &               0 &               2 &               2\\
      MNISTS5LS & 494 &               2 &               2 &               0 &               2 &               2\\
       MODBEALE & 200 &             --- &             --- &             --- &             --- &             ---\\
         MOREBV & 1000 &             333 &            1211 &               0 &             333 &            1211\\
       MSQRTALS & 529 &            3855 &            7337 &               0 &            3855 &            7337\\
       MSQRTBLS & 529 &            2671 &            5099 &               0 &            2671 &            5099\\
          NCB20 & 1010 &             929 &            4353 &               0 &             929 &            4353\\
         NCB20B & 1000 &             --- &             --- &             --- &             --- &             ---\\
       NONCVXU2 & 1000 &            4819 &           24493 &             714 &            5246 &           27413\\
       NONCVXUN & 1000 &            4014 &           20743 &             516 &            5863 &           30231\\
         NONDIA & 1000 &               5 &               5 &               2 &               5 &               5\\
       NONDQUAR & 1000 &             638 &            3400 &             108 &             741 &            3936\\
       NONMSQRT & 529 &             --- &             --- &             --- &             --- &             ---\\
       OSBORNEB & 11 &             232 &             420 &               0 &             232 &             420\\
       OSCIGRAD & 1000 &              60 &              97 &               0 &              60 &              97\\
       OSCIPATH & 100 &             --- &             --- &             --- &             --- &             ---\\
         PARKCH & 15 &            1205 &            3844 &               0 &            1205 &            3844\\
       PENALTY1 & 1000 &               6 &              25 &               4 &               6 &              25\\
       PENALTY2 & 200 &             --- &             --- &             --- &             --- &             ---\\
       PENALTY3 & 200 &              75 &             190 &               1 &              75 &             163\\
       POWELLSG & 1000 &              30 &              70 &              24 &              33 &             122\\
          POWER & 1000 &             159 &             806 &              34 &             586 &            4460\\
         QUARTC & 1000 &             311 &            1838 &               6 &             511 &            2973\\
        SBRYBND & 1000 &             --- &             --- &             --- &             --- &             ---\\
       SCHMVETT & 1000 &             219 &            1121 &               0 &             219 &            1121\\
        SCOSINE & 1000 &             --- &             --- &             --- &             --- &             ---\\
       SCURLY10 & 1000 &             124 &             675 &              13 &             544 &            3986\\
       SCURLY20 & 1000 &             122 &             720 &              17 &             800 &            6386\\
       SCURLY30 & 1000 &              49 &             306 &               7 &             154 &            1088\\
        SENSORS & 100 &              62 &             151 &               0 &              62 &             151\\
        SINQUAD & 1000 &             --- &             --- &             --- &              23 &              78\\
       SPARSINE & 1000 &             474 &             474 &               0 &             474 &             474\\
       SPARSQUR & 1000 &              67 &             294 &               4 &             111 &             580\\
       SPMSRTLS & 1000 &             167 &             369 &               0 &             167 &             369\\
       SROSENBR & 1000 &              12 &              35 &               9 &              12 &              35\\
       SSBRYBND & 1000 &            4932 &           28846 &               0 &            4932 &           28846\\
       SSCOSINE & 1000 &             --- &             --- &             --- &             --- &             ---\\
        STRATEC & 10 &             --- &             --- &             --- &             --- &             ---\\
       TESTQUAD & 1000 &            2499 &           12256 &              12 &            1725 &            7207\\
       TOINTGOR & 50 &             237 &             254 &               0 &             237 &             254\\
       TOINTGSS & 1000 &              17 &              52 &               1 &              17 &              52\\
       TOINTPSP & 50 &             108 &             186 &               0 &             108 &             186\\
       TOINTQOR & 50 &              73 &             100 &               0 &              73 &             100\\
       TQUARTIC & 1000 &              14 &              53 &              10 &              15 &              56\\
         TRIDIA & 1000 &             344 &             354 &               0 &             344 &             354\\
         VARDIM & 200 &               2 &               2 &               0 &               2 &               2\\
       VAREIGVL & 100 &              34 &              52 &               0 &              34 &              52\\
         WATSON & 31 &             594 &            2447 &               0 &             594 &            2447\\
          WOODS & 1000 &              20 &              58 &               1 &              20 &              58\\
        YATP1LS & 120 &              47 &             203 &              43 &              38 &             188\\
        YATP2LS & 120 &              10 &              19 &               7 &              11 &              20\\

    \hline
  \etabular
  }
\etable

\begin{figure}[ht]
  \centering
  \begin{subfigure}[b]{0.48\textwidth}
    \centering
    \includegraphics[width=0.7\textwidth,clip=true,trim=20 0 40 15]{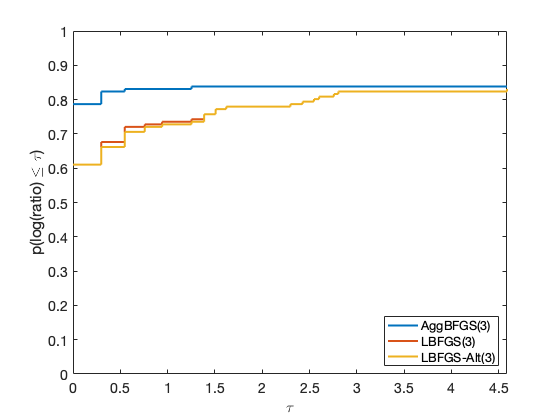}
    \caption{Dolan and Mor\'e Profile, $m = 3$}
  \end{subfigure}
  ~
  \begin{subfigure}[b]{0.48\textwidth}
    \centering
    \includegraphics[width=0.7\textwidth,clip=true,trim=20 0 40 15]{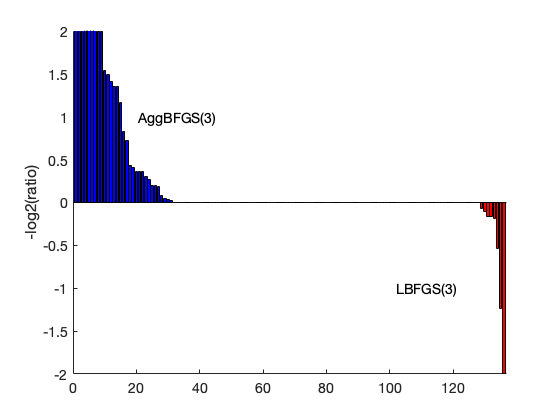}
    \caption{Morales Outperformance Profile, $m = 3$}
  \end{subfigure} \\
  \medskip
  
  \begin{subfigure}[b]{0.48\textwidth}
    \centering
    \includegraphics[width=0.7\textwidth,clip=true,trim=20 0 40 15]{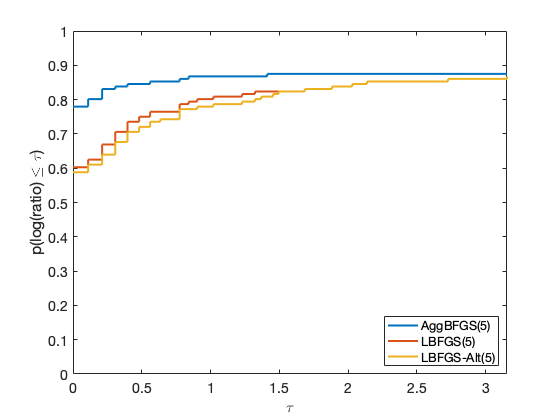}
    \caption{Dolan and Mor\'e Profile, $m = 5$}
  \end{subfigure}
  ~
  \begin{subfigure}[b]{0.48\textwidth}
    \centering
    \includegraphics[width=0.7\textwidth,clip=true,trim=20 0 40 15]{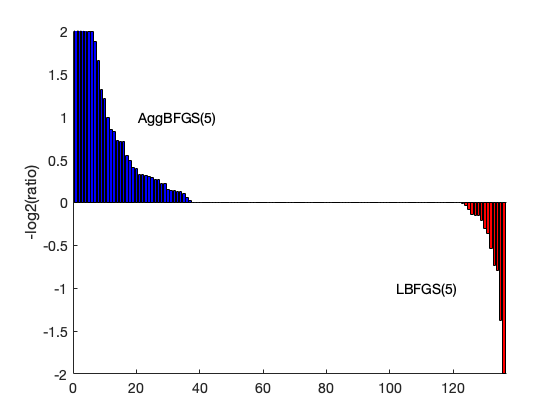}
    \caption{Morales Outperformance Profile, $m = 5$}
  \end{subfigure} \\
  \medskip
  
  \begin{subfigure}[b]{0.48\textwidth}
    \centering
    \includegraphics[width=0.7\textwidth,clip=true,trim=20 0 40 15]{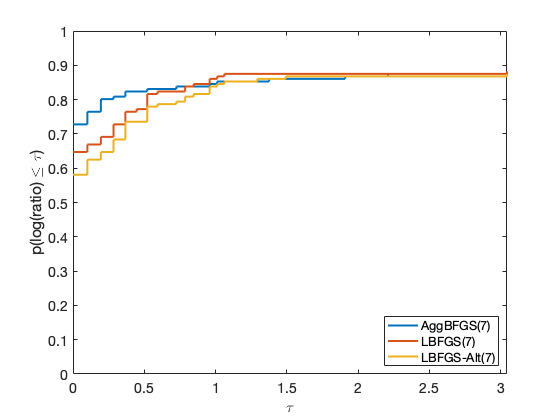}
    \caption{Dolan and Mor\'e Profile, $m = 7$}
  \end{subfigure}
  ~
  \begin{subfigure}[b]{0.48\textwidth}
    \centering
    \includegraphics[width=0.7\textwidth,clip=true,trim=20 0 40 15]{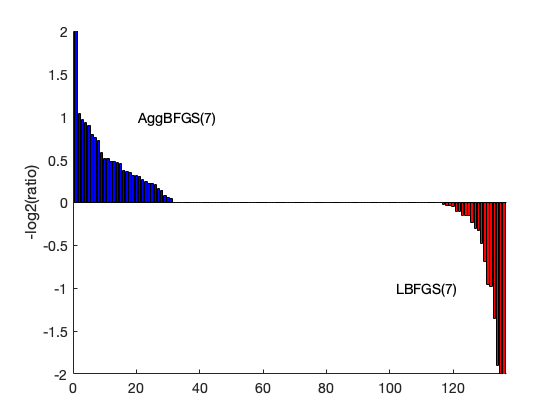}
    \caption{Morales Outperformance Profile, $m = 7$}
  \end{subfigure} \\
  \medskip
  
  \begin{subfigure}[b]{0.48\textwidth}
    \centering
    \includegraphics[width=0.7\textwidth,clip=true,trim=20 0 40 15]{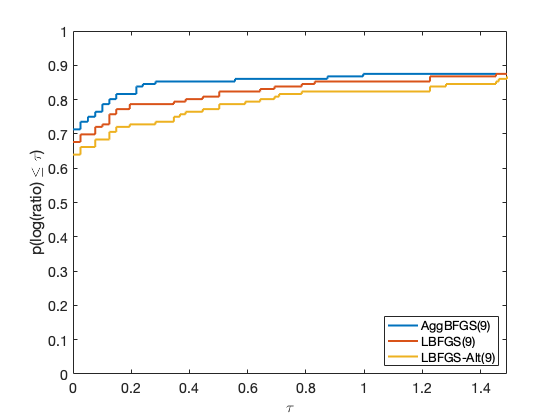}
    \caption{Dolan and Mor\'e Profile, $m = 9$}
  \end{subfigure}
  ~
  \begin{subfigure}[b]{0.48\textwidth}
    \centering
    \includegraphics[width=0.7\textwidth,clip=true,trim=20 0 40 15]{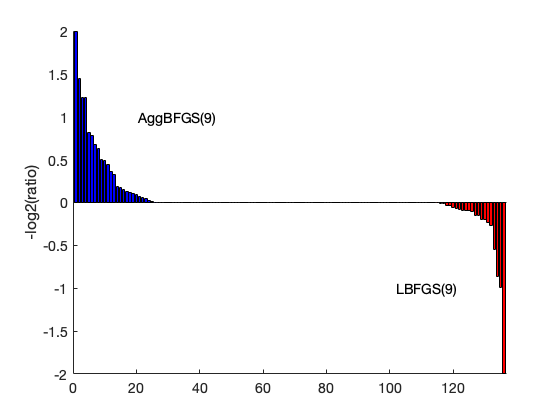}
    \caption{Morales Outperformance Profile, $m = 9$}
  \end{subfigure}
  \caption{Performance profiles for iterations.}
  \label{fig.iterations}
\end{figure}

\begin{figure}[ht]
  \centering
  \begin{subfigure}[b]{0.48\textwidth}
    \centering
    \includegraphics[width=0.7\textwidth,clip=true,trim=20 0 40 15]{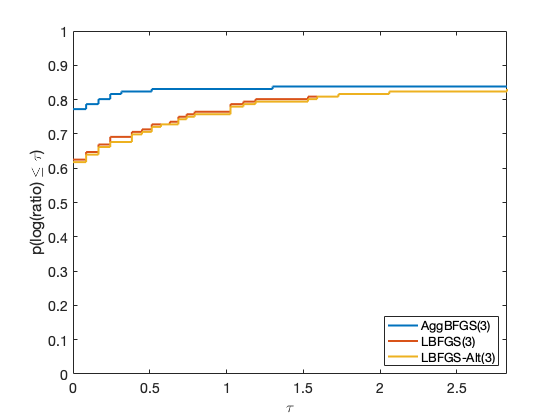}
    \caption{Dolan and Mor\'e Profile, $m = 3$}
  \end{subfigure}
  ~
  \begin{subfigure}[b]{0.48\textwidth}
    \centering
    \includegraphics[width=0.7\textwidth,clip=true,trim=20 0 40 15]{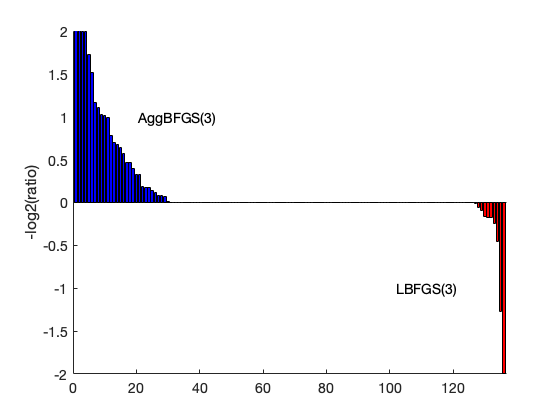}
    \caption{Morales Outperformance Profile, $m = 3$}
  \end{subfigure} \\
  \medskip
  
  \begin{subfigure}[b]{0.48\textwidth}
    \centering
    \includegraphics[width=0.7\textwidth,clip=true,trim=20 0 40 15]{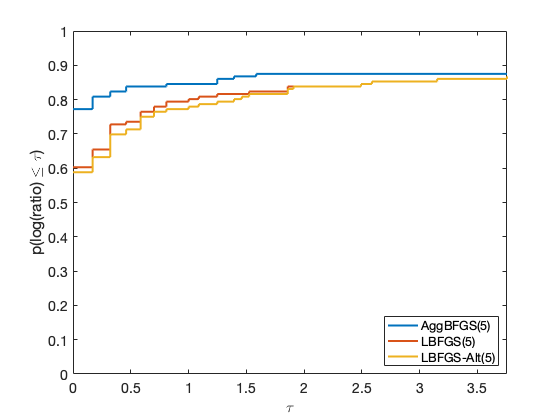}
    \caption{Dolan and Mor\'e Profile, $m = 5$}
  \end{subfigure}
  ~
  \begin{subfigure}[b]{0.48\textwidth}
    \centering
    \includegraphics[width=0.7\textwidth,clip=true,trim=20 0 40 15]{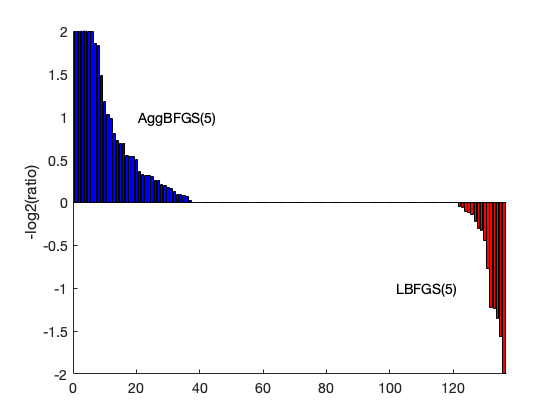}
    \caption{Morales Outperformance Profile, $m = 5$}
  \end{subfigure} \\
  \medskip
  
  \begin{subfigure}[b]{0.48\textwidth}
    \centering
    \includegraphics[width=0.7\textwidth,clip=true,trim=20 0 40 15]{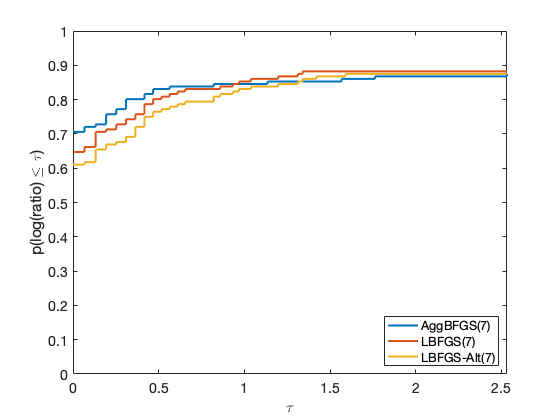}
    \caption{Dolan and Mor\'e Profile, $m = 7$}
  \end{subfigure}
  ~
  \begin{subfigure}[b]{0.48\textwidth}
    \centering
    \includegraphics[width=0.7\textwidth,clip=true,trim=20 0 40 15]{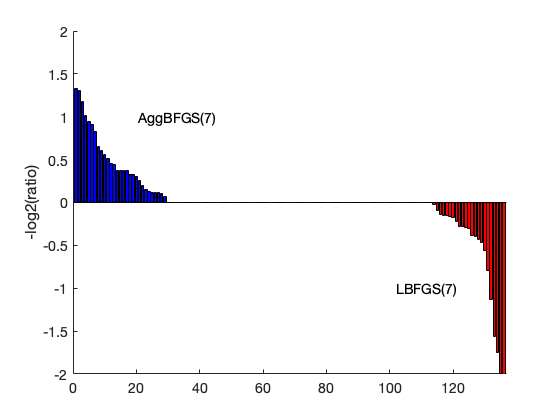}
    \caption{Morales Outperformance Profile, $m = 7$}
  \end{subfigure} \\
  \medskip
  
  \begin{subfigure}[b]{0.48\textwidth}
    \centering
    \includegraphics[width=0.7\textwidth,clip=true,trim=20 0 40 15]{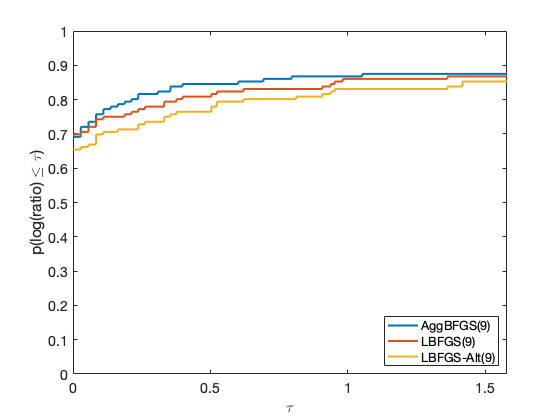}
    \caption{Dolan and Mor\'e Profile, $m = 9$}
  \end{subfigure}
  ~
  \begin{subfigure}[b]{0.48\textwidth}
    \centering
    \includegraphics[width=0.7\textwidth,clip=true,trim=20 0 40 15]{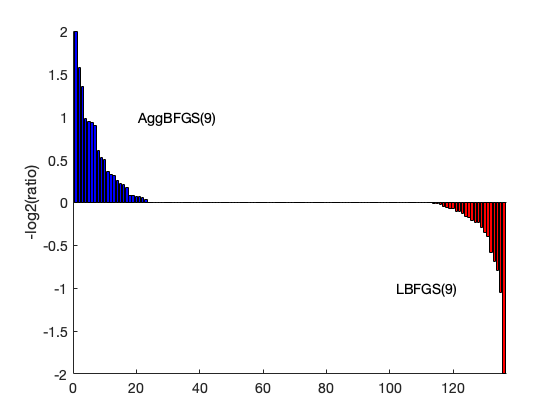}
    \caption{Morales Outperformance Profile, $m = 9$}
  \end{subfigure}
  \caption{Performance profiles for function evaluations.}
  \label{fig.functions}
\end{figure}

\rev{One finds in our results that, generally speaking, \texttt{AggBFGS($m$)} outperforms \texttt{LBFGS($m$)}.  By contrast, the third algorithm in our experiments, which we refer to as \texttt{LBFGS-Alt($m$)}, if anything often performs worse than even \texttt{LBFGS($m$)}.}  This can be seen more clearly in performance profiles.  In Figures~\ref{fig.iterations} and~\ref{fig.functions}, we present the results in the form of two types of profiles: Dolan and Mor\'e performance profiles \cite{DolaMore02} and Morales outperforming factor profiles \cite{Mora02}.  The former type of profile has a graph for each algorithm, where if a graph for an algorithm passes through the point $(\alpha,0.\beta)$, where $\beta$ is a two-digit integer, then on $\beta$\% of the problems the measure required by the algorithm was less than $2^\alpha$ times the measure required by the best algorithm in terms of the measure.  In this manner, an algorithm has performed better than the other if its graph is above and to the left of the graph of the other algorithm.  The second type of profile shows a bar plot of logarithmic outperforming factors, with a bar for each problem in the test set.  In our case, \rev{ignoring the performance of the third algorithm (\texttt{LBFGS-Alt($m$)})}, the bar for each problem takes on the value $-\log_2(meas_{\texttt{AggBFGS($m$)}}/meas_{\texttt{LBFGS($m$)}})$, where $meas_{\texttt{AggBFGS($m$)}}$ is the performance measure for \texttt{AggBFGS($m$)} and $meas_{\texttt{LBFGS($m$)}}$ is the performance measure for \texttt{LBFGS($m$)}.  An upward-pointing bar shows by how much \texttt{AggBFGS($m$)} outperformed \texttt{LBFGS($m$)} on a particular problem, and vice versa for the downward-pointing bars.  We sort the bars \rev{by} value for ease of visualization.  Figure~\ref{fig.iterations} shows profiles using the number of iterations as the performance measure, and Figure~\ref{fig.functions} shows profiles using the number of function evaluations required as the performance measure.  \rev{One may conclude from the profiles that \AggBFGS{} helps performance more often than not, with the benefits being even more significant for smaller values of $m$.}

Much remains to be investigated in terms of the practical use of \AggBFGS{}, including perhaps more sophisticated techniques for implementing the scheme, handling numerical errors, and tuning parameters so that the results may be even better on a wider variety of problems.  Our preliminary experiments in this section motivate such further investigations into the practical use of \AggBFGS.

\section{Conclusion}\label{sec.conclusion}

We have presented a technique for aggregating the curvature pair information in a limited-memory-type BFGS approach such that the corresponding Hessian (and inverse Hessian) approximations are the same as those that would be computed in a full-memory BFGS approach.  The key idea is that if one finds that a stored iterate displacement vector lies in the span of subsequent iterate displacements, then the gradient displacement vectors can be modified in such a way that the pair involving the linearly dependent iterate displacement can be removed with no information being lost.  To the best of our knowledge, this is the first limited-memory-type approach that can behave equivalently to a full-memory method, meaning that it can offer all theoretical properties of a full-memory scheme, such as local superlinear guarantees for the (outer) optimization method.  We have also shown that the application of our aggregation scheme within an L-BFGS method can lead to performance gains over standard L-BFGS.

Our methodology could be extended to other quasi-Newton schemes.  In particular, by the well-known symmetry between the BFGS and DFP updating, it is straightforward to extend our approach for DFP.  In particular, with DFP, rather than looking for linear dependence between iterate displacements, one should look for linear dependence between gradient displacements.  If linear dependence is observed, then one can aggregate the iterate displacements to obtain
\bequation\label{eq.Stilde_formula}
  \Stilde_{1:m} = M^{-1}Y_{1:m} \bbmatrix A & 0 \ebmatrix + s_0 \bbmatrix b \\ 0 \ebmatrix^T + S_{1:m}
\eequation
(cf.~\eqref{eq.Ytilde_formula}) such that $\DFP(M,S_{0:m},Y_{0:m}) = \DFP(M,\Stilde_{1:m},Y_{1:m})$ for $M \succ 0$.  There might also be opportunities for extending our approach for the other members of the Broyden class of updates, although such extensions are not as straightforward.  Indeed, for members of the class besides BFGS and DFP, one likely needs to modify both iterate and gradient displacements to aggregate curvature information.

There are also other opportunities for designing practical adaptations of our scheme.  Perhaps the most straightforward idea is the following: Suppose that one is employing an L-BFGS$(m)$ approach, one has $m$ curvature pairs already stored, and one performs a new optimization algorithm iteration to yield a new curvature pair.  Rather than simply discard the oldest stored curvature pair if a previous iterate displacement does not lie in the span of subsequent displacements, one could project this pair's iterate displacement onto the span of the subsequent displacements (and possibly project the pair's gradient displacement onto a subspace), then apply our aggregation scheme.  This offers the opportunity to maintain more historical curvature information while still only storing/employing at most $m$ pairs of vectors.  One might also imagine other approaches for projecting information into smaller-dimensional subspaces in order to employ our scheme, rather than simply discarding old information.  Such techniques might not attain the theoretical properties of a full-memory approach, but could lead to practical benefits.

One could employ our aggregation scheme with no modifications necessary if one employs an optimization algorithm that intentionally computes sequences of steps in low-dimensional subspaces, such as in block-coordinate descent.  Another setting of interest is large-scale constrained optimization where the number of degrees of freedom (i.e., number of variables minus the number of active constraints) is small relative to $n$.  For such a problem, various algorithms compute search directions that lie in subspaces that are low-dimensional relative to $n$.

\bibliographystyle{plain}
\bibliography{references}

\begin{thebibliography}{10}

\bibitem{berahas2016multi}
A.~S. Berahas, J.~Nocedal, and M.~Tak{\'a}{\v{c}}.
\newblock {A multi-batch L-BFGS method for machine learning}.
\newblock In {\em Advances in Neural Information Processing Systems}, pages
  1055--1063, 2016.

\bibitem{berahas2020robust}
A.~S. Berahas and M.~Tak{\'a}{\v{c}}.
\newblock A robust multi-batch l-bfgs method for machine learning.
\newblock {\em Optimization Methods and Software}, 35(1):191--219, 2020.

\bibitem{BoggByrd19}
P.~T. Boggs and R.~H. Byrd.
\newblock Adaptive, limited-memory bfgs algorithms for unconstrained
  optimization.
\newblock {\em SIAM Journal on Optimization}, 29(2):1282--1299, 2019.

\bibitem{BonnGilbLemaSaga95}
J.~F. Bonnans, J.~Ch. Gilbert, C.~Lemar{\'e}chal, and C.~A. Sagastiz{\'a}bal.
\newblock {A family of variable metric proximal methods}.
\newblock {\em Mathematical Programming}, 68(1):15--47, 1995.

\bibitem{Broy70}
C.~G. Broyden.
\newblock {The convergence of a class of double-rank minimization algorithms}.
\newblock {\em Journal of the Institute of Mathematics and Its Applications},
  6(1):76--90, 1970.

\bibitem{byrd2016stochastic}
R.~H. Byrd, S.~L. Hansen, J.~Nocedal, and Y.~Singer.
\newblock {A stochastic quasi-Newton method for large-scale optimization}.
\newblock {\em SIAM Journal on Optimization}, 26(2):1008--1031, 2016.

\bibitem{ByrdNoce89}
R.~H. Byrd and J.~Nocedal.
\newblock {A tool for the analysis of quasi-Newton methods with application to
  unconstrained minimization}.
\newblock {\em SIAM Journal on Numerical Analysis}, 26(3):727--739, 1989.

\bibitem{ByrdNoceSchn94}
R.~H. Byrd, J.~Nocedal, and R.~B. Schnabel.
\newblock {Representations of quasi-Newton matrices and their use in limited
  memory methods}.
\newblock {\em Mathematical Programming}, 63:129--156, 1994.

\bibitem{ByrdNoceYuan87}
R.~H. Byrd, J.~Nocedal, and Y.~Yuan.
\newblock {Global convergence of a class of quasi-Newton methods on convex
  problems}.
\newblock {\em SIAM Journal on Numerical Analysis}, 24(5):1171--1189, 1987.

\bibitem{Curt16}
F.~E. Curtis.
\newblock {A self-correcting variable-metric algorithm for stochastic
  optimization}.
\newblock In {\em Proceedings of the 48th International Conference on Machine
  Learning}, volume~48, pages 632--641, New York, New York, USA, 2016.
  Proceedings of Machine Learning Research.

\bibitem{CurtQue13}
F.~E. Curtis and X.~Que.
\newblock {An adaptive gradient sampling algorithm for nonsmooth optimization}.
\newblock {\em Optimization Methods and Software}, 28(6):1302--1324, 2013.

\bibitem{CurtQue15}
F.~E. Curtis and X.~Que.
\newblock {A quasi-Newton algorithm for nonconvex, nonsmooth optimization with
  global convergence guarantees}.
\newblock {\em Mathematical Programming Computation}, 7:399--428, 2015.

\bibitem{CurtRobiZhou19}
F.~E. Curtis, D.~P. Robinson, and B.~Zhou.
\newblock {A self-correcting variable-metric algorithm framework for nonsmooth
  optimization}.
\newblock {\em IMA Journal of Numerical Analysis}, to appear, 2019.

\bibitem{Davi91}
W.~C. Davidon.
\newblock {Variable metric method for minimization}.
\newblock {\em SIAM Journal on Optimization}, 1(1):1--17, 1991.

\bibitem{DennMore74}
J.~E. Dennis and J.~J. Mor\'e.
\newblock {A characterization of superlinear convergence and its application to
  quasi-Newton methods}.
\newblock {\em Mathematics of Computation}, 28(126):549--560, 1974.

\bibitem{DennSchn96}
J.~E. Dennis and R.~B. Schnabel.
\newblock {\em {Numerical Methods for Unconstrained Optimization and Nonlinear
  Equations}}.
\newblock Society for Industrial and Applied Mathematics (SIAM), Philadelphia,
  PA, USA, 1996.

\bibitem{DolaMore02}
E.~D. Dolan and J.~J. More.
\newblock {Benchmarking Optimization Software with Performance Profiles}.
\newblock {\em Mathematical Programming}, 91(2):201--213, 2002.

\bibitem{Flet70}
R.~Fletcher.
\newblock {A new approach to variable metric algorithms}.
\newblock {\em Computer Journal}, 13(3):317--322, 1970.

\bibitem{GillGoluMurrSaun74}
P.~E. Gill, G.~H. Golub, W.~Murray, and M.~A. Saunders.
\newblock {Methods for modifying matrix factorizations}.
\newblock {\em Mathematics of Computation}, 126(28):505--535, 1974.

\bibitem{Gold70}
D.~Goldfarb.
\newblock {A family of variable metric updates derived by variational means}.
\newblock {\em Mathematics of Computation}, 24(109):23--26, 1970.

\bibitem{GoulOrbaToin15}
N.~I.~M. Gould, D.~Orban, and Ph.~L. Toint.
\newblock {CUTEst: A constrained and unconstrained testing environment with
  safe threads for mathematical optimization}.
\newblock {\em Computational Optimization and Applications}, 60(3):545--557,
  2015.

\bibitem{gower2016stochastic}
R.~Gower, D.~Goldfarb, and P.~Richt{\'a}rik.
\newblock {Stochastic block BFGS: Squeezing more curvature out of data}.
\newblock In {\em International Conference on Machine Learning}, pages
  1869--1878, 2016.

\bibitem{HaarMietMaek04}
N.~Haarala, K.~Miettinen, and M.~M. M{\"a}kel{\"a}.
\newblock {New limited memory bundle method for large-scale nonsmooth
  optimization}.
\newblock {\em Optimization Methods and Software}, 19(6):673--692, 2004.

\bibitem{keskar2016adaqn}
N.~S. Keskar and A.~S. Berahas.
\newblock {ADAQN: An adaptive quasi-Newton algorithm for training RNNs}.
\newblock In {\em Joint European Conference on Machine Learning and Knowledge
  Discovery in Databases}, pages 1--16. Springer, 2016.

\bibitem{KoldOLeaNaza98}
T.~G. Kolda, D.~P. O'Leary, and L.~Nazareth.
\newblock {BFGS with Update Skipping and Varying Memory}.
\newblock {\em SIAM Journal on Optimization}, 8(4):1060--1083, 1998.

\bibitem{LewiOver13}
A.~S. Lewis and M.~L. Overton.
\newblock {Nonsmooth optimization via quasi-Newton methods}.
\newblock {\em Mathematical Programming}, 141(1):135--163, 2013.

\bibitem{MiffSunQi98}
R.~Mifflin, D.~Sun, and L.~Qi.
\newblock {Quasi-Newton bundle-type methods for nondifferentiable convex
  optimization}.
\newblock {\em SIAM Journal on Optimization}, 8(2):583--603, 1998.

\bibitem{mokhtari2015global}
A.~Mokhtari and A.~Ribeiro.
\newblock {Global convergence of online limited memory BFGS}.
\newblock {\em The Journal of Machine Learning Research}, 16(1):3151--3181,
  2015.

\bibitem{Mora02}
J.~L. Morales.
\newblock {A Numerical Study of Limited Memory BFGS Methods}.
\newblock {\em Applied Mathematics Letters}, 15:481--487, 2002.

\bibitem{Noce80}
J.~Nocedal.
\newblock {Updating quasi-Newton matrices With limited storage}.
\newblock {\em Mathematics of Computation}, 35(151):773--782, 1980.

\bibitem{NoceWrig06}
J.~Nocedal and S.~J. Wright.
\newblock {\em {Numerical Optimization}}.
\newblock Springer New York, {Second} edition, 2006.

\bibitem{Pear69}
J.~D. Pearson.
\newblock {Variable metric methods of minimisation}.
\newblock {\em The Computer Journal}, 12(2):171--178, 1969.

\bibitem{Powe76}
M.~J.~D. Powell.
\newblock {Some global convergence properties of a variable metric algorithm
  for minimization with exact line searches}.
\newblock In R.~W. Cottle and C.~E. Lemke, editors, {\em Nonlinear Programming,
  SIAM-AMS Proceedings, Vol.~IX}, Harwell, England, 1976. American Mathematical
  Society.

\bibitem{Ritt79}
K.~Ritter.
\newblock {Local and superlinear convergence of a class of variable metric
  methods}.
\newblock {\em Computing}, 23(3):287--297, 1979.

\bibitem{Ritt81}
K.~Ritter.
\newblock {\em {Global and superlinear convergence of a class of variable
  metric methods}}, pages 178--205.
\newblock Springer Berlin Heidelberg, Berlin, Heidelberg, 1981.

\bibitem{Rose60}
H.~H. Rosenbrock.
\newblock {An automatic method for finding the greatest or least value of a
  function}.
\newblock {\em The Computer Journal}, 3(3):175--184, 1960.

\bibitem{schraudolph2007stochastic}
N.~N. Schraudolph, J.~Yu, and S.~G{\"u}nter.
\newblock {A stochastic quasi-Newton method for online convex optimization}.
\newblock In {\em Artificial Intelligence and Statistics}, pages 436--443,
  2007.

\bibitem{Shan70}
D.~F. Shanno.
\newblock {Conditioning of quasi-Newton methods for function minimization}.
\newblock {\em Mathematics of Computation}, 24(111):647--656, 1970.

\bibitem{VlceLuks01}
J.~Vl{\v{c}}ek and L.~Luk{\v{s}}an.
\newblock {Globally convergent variable metric method for nonconvex
  nondifferentiable unconstrained minimization}.
\newblock {\em Journal of Optimization Theory and Applications},
  111(2):407--430, 2001.

\bibitem{wang2017stochastic}
X.~Wang, S.~Ma, D.~Goldfarb, and W.~Liu.
\newblock {Stochastic quasi-Newton methods for nonconvex stochastic
  optimization}.
\newblock {\em SIAM Journal on Optimization}, 27(2):927--956, 2017.

\end{thebibliography}

\end{document}